\documentclass[10pt]{amsart}
\usepackage{amsmath}
 \newcommand{\be}{\begin{enumerate}}
\newcommand{\ee}{\end{enumerate}} \newcommand{\bd}{\begin{description}}
 \newcommand{\ed}{\end{description}} \newcommand{\bi}{\begin{itemize}}
\newcommand{\ei}{\end{itemize}}
\newtheorem{thm}{Theorem}[section] %

\newtheorem{corollary}[thm]{Corollary} %
\numberwithin{equation}{section}
\newtheorem{lemma}[thm]{Lemma} %

\newcommand{\eqn}{\begin{eqnarray}}
\newcommand{\eeqn}{\end{eqnarray}}

\begin{document}

\title{Enumerative Combinatorics of Simplicial and Cell Complexes: Kirchhoff and Trent Type Theorems}

\author{Sylvain E. Cappell and Edward Y. Miller}

\maketitle

\centerline{ }

\vspace{-.4in}

\centerline{ABSTRACT}

{ \small
This paper considers three separate matrices associated to graphs
and  (each dimension of) cell complexes. It relates all the coefficients of their respective
characteristic polynomials to the geometric and combinatorial enumeration
of three kinds of subobjects. The matrices are:
the mesh matrix for  integral d-cycles, the mesh matrix for integral d-boundaries, and
the Kirchhoff matrix, i.e., the combinatorial Laplacian, for integral (d-1)-chains.

Trent's theorem states that  the determinant of
the mesh matrix on 1-cycles of  a connected graph is equal to the number of spanning trees \cite{Trent,Trent1}.
 Here this theorem is extended
to the mesh matrix on $d$ cycles in an arbitrary cell complex and, new even for graphs,
to enumerative combinatorial interpretation of all of  the coefficients of  its characteristic
polynomial.
This last  is well defined once a basis for the integral $d$-cycles
is chosen.
Additionally, a parallel result for  the mesh matrix for integral d-boundaries  is proved.

Kirchhoff's theorem states that the product of the non-zero eigenvalues
of the Kirchhoff matrix, i.e., combinatorial Laplacian, for connected graphs equals
n times the number of spanning trees with n the number of vertices.
 Lyons has generalized Kirchhoff's result on the product of the
non-zero eigenvalues of the Kirchhoff or combinatorial Laplacian on (d-1)-chains
to cell complexes for $d >1$ \cite{Lyons}. The present analysis extends this to all coefficients
of the characteristic polynomial.

An evaluation of the Reidemeister-Franz torsion of the cell complex
with respect to its integral basis gives relations between these combinatorial
invariants.

\vspace{.1in}

\begin{tabular}{ll}
 Section 1. & Introduction \\

 Section 2. & Theorem on mesh matrix for $d$ cycles: evaluation of its  characteristic polynomial \\

Section 3.  & Proof of the main theorem for mesh matrix for $d$ cycles: theorem \ref{thm1} \\

Section 4.  &  Alternative version of the squared covolume $det(j_V \cdot j_V^\star)$, cases $ k \ge 0$ \\

Section 5. &  The mesh matrix for $d$ boundaries: evaluation of its  characteristic polynomial \\

Section 6.  & Proof of  main theorem \ref{thm2} for mesh matrix for $d$ boundaries: theorem  \\

Section 7: & Evaluation of all coefficients of mesh matrices for cycles and boundaries \\ &
 with respect to   geometrically defined bases. \\

Section 8.  &  A Kirchhoff type  theorem for cell  complexes ( extending    Lyons \cite{Lyons} )
\\

Section 9.  & Homology covolume squared via determinants of mesh matrices  \\

Section 10. & Combinatorial identities via a computation of Reidemeister-Franz torsion \\

Section 11. & Weights, Variant Settings, and Examples \\

Section 12. & Enumerative combinatorics of the simplex   ( extending    Kalai \cite{Kalai} )  \\

\end{tabular}

}

{ \Large
\section{\textbf{Introduction}}

\vspace{.1in}

This paper considers three separate matrices associated to graphs
and  (each dimension of) cell complexes. It relates all the coefficients of their respective
characteristic polynomials to the geometric and combinatorial enumeration
of three kinds of sub-objects.
Their determinants and more generally their  characteristic polynomials are interpreted
in terms of correspondingly three different aspects of the enumerative combinatorics of the graph or cell complex
(in each dimension).

The oldest such matrix,  of the three studied here, is Kirchhoff's
matrix, which is the combinatorial Laplacian $\Delta_0$ \cite{Kirchhoff,Chung}
on 0-chains. For graphs, $\Delta_0$ is the degree matrix minus the adjacency
matrix.
For a connected graph with $n$ vertices, Kirchhoff showed that the product
of the non-vanishing eigenvalues of $\Delta_0$ is equal to $n$
times the number of spanning trees.  This
interpretation of the product of the non-vanishing eigenvalues
of the combinatorial Laplacian on (d-1)-chains  was generalized to higher dimensional
complexes by Lyons \cite{Lyons} as an explicit weighted sum over
pairs of combinatorial objects which he called
bases and cobases, and here (consistent with the classical terminology for graphs) are  called (d-dimensional) spanning forests and ((d-1)-dimensional) spanning coforests.
 In the restricted context of
 d-dimensional complexes which are rationally acyclic in dimensions less than d,
 parallel theorems to those of Lyons on the product of the eigenvalues  of the Krichhoff matrix
  had earlier appeared in the work of  Duval, Klivans, and Martin  \cite{DuvalKlivansMartin}.

In the case of graphs the complete characteristic polynomial of the Kirchhoff matrix
had been studied and each of  its coefficients had been evaluated by an   enumerative
combinatorial formula, a sum over special subgraphs, forests with
specified number of edges with associated weights \cite{ChaikenKleitman,Chaiken}. For  cell complexes,
analysis of the characteristic polynomials is  developed here giving
an analogous interpretation for all coefficients  of the polynomial in each dimension.

The second type of matrix studied here is the  higher dimensional analogue of the
mesh matrix introduced by Trent \cite{Trent,Trent1} to study enumerative combinatorics
of graphs. It is based on the idea that the integral cycles
have a natural pairing to the integers,  yielding a mesh matrix once
a basis is chosen. As shown by Trent, the determinant of the mesh matrix of a connected
graph equals the number of spanning trees. Here  this classical work is extended
in two new directions. One extends the definition and results about the
Trent matrix for graphs to simplicial or cell complexes, obtaining
a mesh matrix in each dimension $d$ with determinant  a sum over (d-dimensional) spanning forests with
explicit weights. Next, we  give a   geometrical interpretation
for each coefficient  in the characteristic polynomial of the mesh matrix
of graphs as well as in the higher dimensional context. The mesh matrix is defined
for each choice of integral basis for the d-cycles;  its determinant
is independent of this basis, while the  coefficients of other powers of $t$
in  the characteristic polynomial depend
on this choice.  We see that in \S 5, the coefficient of each power of $t$  is shown to be  sum   over geometrical
objects, (d-dimensional) k augmented  spanning forests, with positive integer weights depending
in a specified way on the
chosen basis. This is of interest and novel even in  the graph setting. Moreover, these weights are here
given for cell complexes and graphs several topological and combinatorial interpretations
in  settings of  geometrical interest.

Thirdly, the integral boundaries of chains in each dimension also have a natural
pairing to the integers. Once an integral basis is
chosen a matrix arises. Its determinant is
the sum over (d-dimensional) spanning coforests with  explicit positive integral  weights
which are studied below.
Again each coefficient of powers of $t$ in the characteristic polynomial
has an expression as a sum uniformly over (d-dimensional)  k reduced coforests
with certain positive integral weights.

The geometrical objects being enumerated with weights in the
formulae for each of the three  settings are thus  different.
Let $n_d = dim\ C_d(X;R), z_d = dim \ Z_d(X;R)$, $b_d = dim \ B_d(X;R)$,
the dimensions of the $d$-chains, $d$-cycles, $d$-boundaries.
For the reader's convenience, these are tabulated in the following
table.

}

{\small
\vspace{.0in}
\hspace{.2in}
\begin{tabular}{c||c|c|c|c}

& \multicolumn{2}{|c|}{for   connected graphs \ (d=1)  } & \multicolumn{2}{|c}{for  cell  complexes } \\
\hline \hline

Matrix & product & coefficient of $t^L$   & product & coefficient of $t^{L}$ \\
Types & non-vanishing & of characteristic   & non-vanishing &  of characteristic \\
& eigenvalues  &  polynomial   & eigenvalues  & polynomial \\
 &  & & $ $  & \\
 \hline \hline

 &   \multicolumn{4}{|c}{ \textbf{ Geometrical objects summed over to yield algebraic invariant}} \\
 \hline \hline

Kirchhoff's  & spanning & forests of size m &  pairs, spanning  & pairs, forest   \\
Combinatorial &               trees  &     $ L=n_0-m$             &    forest and    & of size m  and   \\
Laplacian            &      $ L=1 $                  &                  &  spanning         & spanning  \\
on (d-1)-chains                     &              &                     &  coforest  in  it   &  coforest in it \\
       $n_{d-1}  \times n_{d-1}$               &             &                      & $L=n_{d-1}- b_{d-1}$  & $L=n_{d-1}-m$ \\
                      &&&& $ 0 \le m \le b_{d-1} $ \\
\hline

Trent  & spanning &  k augmented   &  spanning forests & k augmented   \\
Mesh matrix &    trees             &  spanning forests                 & $L=0$ & spanning forests   \\
on d-cycles          &      $L=0$                 & $L=k$         &  &  $L=k$   \\
$z_d \times z_d$ &&&& $0 \le k <  z_d $  \\
\hline

Mesh matrix & -  & -  & spanning  coforests & k reduced   \\
on d-boundaries &                 &                   & $L=0 $ &   spanning coforests  \\
   $b_d \times b_d$  &&& & $L=k$ \\
    &&&& $0 \le k < b_d$ \\
\hline \hline

\end{tabular}

In the Kirchhoff case the pairs, (forest, spanning coforest in it) are in dimensions (d,d-1) respectively;
in the cycle mesh matrix case the k augmented spanning forest  is in dimension d; in the boundary mesh matrix
case the k reduced spanning coforest is in dimension d. A k augmented spanning forest is a subcomplex
for which removing some choice  of k \  d-dimensional cells yields a spanning forest; a k reduced spanning coforest is
a subcomplex for which adding some choice of k \   d-dimensional cells yields a spanning coforest.

}
\vspace{.1in}
\Large{

Additionally, for the d-cycles given a spanning forest there is a geometrically defined
basis for the rational cycle group $Z_d(X;Q)$; again for the d-boundaries given a spanning forest in dimension $d+1$
there is a rationally defined basis for the rational boundary group $B_d(X;Q)$. Consequently, one obtains
two  types of geometrically defined mesh type matrices. In section 7,
in each of these  cases
the complete characteristic polynomial is evaluated combinatorially.

\vspace{.1in}
A spanning tree for a connected graph with $n$ vertices is a choice of $n-1$ edges which together carry
no 1-cycles. More generally, even when  the graph is not necessarily connected, a choice of spanning tree for each connected
component is called a spanning forest. In  this possibly disconnected
 context, Trent's theorem for graphs asserts  that  the determinant of the mesh matrix for cycles
is the number of spanning forests \cite{Trent,Trent1}. Correspondingly, the  generalization
to the complete characteristic polynomial  below is expressed in terms
of spanning forests and its extension, k-augmented spanning forests.

The mesh matrix of $1$ cycles  for a graph and more generally for $d$ cycles of a $CW$ complex $X$ is
easily explained. The $d$ chains, $C_d(X;Z)$ have a basis of oriented  $d$ cells,
and this yields a  combinatorially defined pairing by setting $<\sigma_j, \sigma_k> = 0$
if the cells are different and $+1$ if $j = k$. By linearity this gives a pairing of integral
chains
$$
<.,.> \ : \ C_d(X;Z) \times C_d(X;Z) \rightarrow Z .
$$
By definition the integral $d$ cycles, $Z_d(X;Z)$ are the kernel of the boundary
mapping $\partial : C_d(Z;Z) \rightarrow C_{d-1}(X;Z)$.
Now including the $d$ cycles $Z_d(X;Z) \subset C_d(X;Z)$ the composite
$Z_d(X;Z) \times Z_d(X;Z) \subset C_d(X;Z) \times C_d(X;Z) \rightarrow Z$
gives the pairing $<.,.>$.

Taking a basis for the $d$ cycles, which is a free abelian group, say $\{ z_r\}$,
there is the associated mesh matrix on $d$ cycles
$$
M\{z_r\} =  \{ < z_j, z_k> \} .
$$
An analogous matrix is associated to  the integral  $d$ boundaries, $B_d(X;Z) = image \ of \
\partial :C_{d+1}(X;Z) \rightarrow C_d(X;Z)$ giving a mesh matrix for $d$ boundaries, once
 a basis  is chosen for the integral boundaries.

\vspace{.1in}
In this paper, the determinant
of the mesh matrix on $d$-cycles is shown   to be a  simple positive integer
 weighted sum over what are  called here ``spanning forests'' in a finite cell  complex $X$,
   extending  the classical notion of spanning forests for graphs.
 These spanning forests are subsets of $d$ cells
 maximal with respect to the property of carrying no $d$ cycles.
 This extends Trents's theorem to each dimension, say d,  of a cell  complex. Two distinct combinatorial evaluations
 of the weights are here provided.

 There is a corresponding notion of k augmented spanning forest  which are just  spanning forests
  with $k$  added $d$ cells.
 Having picked a basis for
the integral $d$-cycles, the mesh matrix on $d$-cycles
 is a well-defined  integral valued $z_d \times z_d$ square matrix.
  It is shown that
the $(z_d-k)^{th}$ elementary symmetric function of the eigenvalues of
of the mesh matrix for cycles   is  a sum over corresponding k augmented spanning forest with  positive integral weights.
These weights are given  an explicit dependence on the basis of $d$-cycles chosen.
This extends Trent's theorem  to the complete characteristic polynomial of the mesh matrix for cycles.
The weights are given in two quite different  formats (but equal).
 This is of interest
even in the graph, $d=1$, case. These definitions and results appear  in \S \ref{sectioncyc}.
( The spanning forests are called bases in Lyons \cite{Lyons}. )

\vspace{.1in}
Additionally, the  mesh matrix for $d$ boundaries has a parallel development for cell  complexes; although it is trivial
for graphs. Having chosen a  basis for the integral boundaries, it is a well defined matrix. The determinant
is obtained by a  positively integer weighted sum over spanning coforests.
Having picked a basis for the integral $d$-boundaries,
the $(b_d-k)^{th}$ elementary symmetric function of the eigenvalues of
the mesh matrix for $d$ boundaries is  shown to be a sum over k  reduced spanning coforests
with appropriate positive integral weights depending on the choice of basis. A k reduced spanning coforest
is just a spanning coforest with  $k $ cells removed.
 These definitions and results appear  in \S \ref{sectionbd}. Again   two distinct evaluations
 of  the weights are provided.

In section 7,
for each of two mesh matrix   cases,
the complete characteristic polynomial is evaluated combinatorially in yet another way, using
a geometrically defined basis of the rational $d$ cycles and boundaries.
This employs a fixed  choice of a spanning forest in dimension $d$ and $d+1$
respectively.

\vspace{.1in}
 In the Kirchhoff case and in the mesh cycle and mesh boundary case
  there is a formal analogy. The combinatorial  Laplacian on $d$ chains,
 the mesh matrix
 for $d$-cycles, and  the mesh matrix for $d$ boundaries
    each are represented by matrices of the form, $A \cdot A^t$.
 The Cauchy-Binet theorem applies as an algebraic step in the analysis
 of such $A \cdot A^t$.  The non-zero eigenvalues
 of $A \cdot A^t$ are real,  non-negative numbers, their positive square roots are often
 referred to as the singular values of the matrix $A$.

This  algebraic procedure applied in the combinatorial Laplacian  setting recapitulates
the approach of Lyons who proved that the product of the non-zero eigenvalues of this
Laplacian, acting on (d-1) chains, equals the a sum over (d dimensional) spanning forests
and ( (d-1) dimensional) spanning coforests contained in them with positive integral weights
in  the terminology used here [bases and cobases in his treatment] \cite{Lyons}. Here this result
is extended to the coefficients of the complete characteristic polynomial using
a sum over pairs, a (d dimensional) forest  and a ($(d-1)$ dimensional) spanning coforest in it.
These results  appear in \S 8.

\vspace{.1in}
There is a large and growing literature on the general topic of generalizations, analogues and applications of
Kirchhoff type theorems.
A small sample of relevant combinatorially oriented papers and monographs would include
(1847) G. Kirchhoff \cite{Kirchhoff},
 (1889) A. Cayley \cite{Cayley},  (1935) H. Whitney \cite{Whitney},
 (1954) H. M. Trent \cite{Trent}, (1955) H. M. Trent \cite{Trent1}, (1960) N. I . Biggs \cite{Biggs},
  (1961) A. Nerode and H. Shank \cite{NerodeShank}, (1967) P. R. Bryant \cite{Bryant},
(1971)   W. K. Chen \cite{Chen1}, (1971) H. Shank \cite{Shank},  (1976) S. B. Maurer \cite{Maurer},
(1976) J.W. Moon \cite{Moon}, (1978) S. Chaiken and D. Kleitman \cite{ChaikenKleitman}, (1979) W. Tutte \cite{Tutte},
 (1982) S. A. Chaiken \cite{Chaiken},
(1983) G. Kalai \cite{Kalai}, (1992) R. M. Adin \cite{Adin}, (2000)  R. Kenyon \cite{Kenyon1},
(2009) A. Duval, C. Klivans, and J. Martin \cite{DuvalKlivansMartin}, (2009) R. Lyons \cite{Lyons},
(2009) A. Petersson,
(2011)  R. Kenyon\cite{Kenyon2}, (2013)  M. Catanzaro, V. Chernyak,
and  J. Klein \cite{CatanzaroChernyakKlein2},  the items after 1983  in   higher dimensional settings.

\vspace{.1in}
The determinants of the mesh matrices for $d$ cycles divided by that for $d$ boundaries
is directly related to the determinant of a  corresponding matrix for the homology.
A computation of the Reidemeister-Franz torsion of the $CW$ complex $X$,
yields interesting relations with   the above combinatorial items. These definitions and combinatorial
identities  appear in \S \ref{sectRF}. A use of Reidemeister-Franz torsion
in the combinatorial context also appears in the work of Catanzaro, Chernyak, Klein \cite{CatanzaroChernyakKlein1,CatanzaroChernyakKlein2}.

 In section \S \ref{sectionKalai}
the example of the standard  simplex is considered, using and extending  results of Kalai \cite{Kalai}.

            This paper emphasizes combinatorial identities. We note, however, that as the coefficients in the
formulas produced are generally shown to be positive integers, a wealth of interesting inequalities, often
involving simplified expressions, result as unstated  corollaries.

      In a future paper we will see how the Cauchy-Binet Theorem can be extended  from the
      context of the determinant, which can be regarded as the trace of a one dimensional
      irreducible representation to the context of the trace of an arbitrary irreducible
      representation \cite{Weyl}. By these means, additional combinatorial-geometric  formulae
      can be derived for graphs and cell complexes.

\section{\textbf{Theorem on mesh matrix for $d$ Cycles}} \label{sectioncyc}
\vspace{.1in}

Let $X$ be a finite  CW-complex.  Orient the cells. This data provides
a basis for the integral $d$-chains, $C_d(X;Z)$, a free abelian group.
This basis is called the cellular basis. Let the
boundary map be  $\partial_d: C_d(X;Z)\rightarrow C_{d-1}(X;Z)$, and set the $d$-cycles $Z_d(X;Z) = kernel \ of \
\partial_d$, and $(d-1)$ boundaries $B_{d-1}(X;Z) = image\ of \ \partial_d$, and correspondingly for
real coefficients $R$.

Define a natural integral valued pairing on
$C_d(X;Z)$ by
$$
\begin{array}{l}
<. , .  > \ : \ C_d(X;Z) \times C_d(X;Z) \rightarrow Z \\
 (\ \Sigma_j \ a_j [\sigma_j], \Sigma_k \  b_k [\sigma_k] \ ) \mapsto \Sigma_j \ a_j \cdot b_j
\end{array}
$$
where $S_d =  \{ [\sigma_j]\}$ is the cellular basis in dimension $d$.
That is, the  elements of the cellular basis  are of  length one and orthogonal.

As above,  let $n_d= |S_d|$ be the number of $d$-cells and  $z_d,b_d$ be the ranks of the free abelian groups, $Z_d(X;Z),
 B_d(X;Z)$, of integral $d$ cycles, $d$ boundaries respectively.

 Introduce the mesh pairings on $d$ cycles and $d$ boundaries:
$$
\begin{array}{l}
\mu : Z_d(X;Z) \times Z_d(X;Z) \subset C_d(X;Z) \times C_d(X;Z) \stackrel{<. \ ,\ .>}{\rightarrow} Z \\
\nu : B_d(X;Z) \times B_d(X;Z) \subset C_d(X;Z) \times C_d(X;Z) \stackrel{<.\ , \ .>}{\rightarrow} Z
\end{array}
$$

Let $d \ge 1$ and  fix an ordered  basis, say $\{z_r\ | \ r=1,\cdots,  z_d\}$ for the
integral  $d$ cycles $Z_d(X;Z)$ and define the associated mesh matrix for $d$ cycles $M\{z_r\}$
a $z_d \times z_d$ integral valued matrix by
$$
M\{z_r\}_{i,j} = < z_i,z_j> .
$$
Similarly, fix an ordered basis, say $\{ b_s \ | \ s=1, \cdots, b_d\}$ for the integral
$d$ boundaries $B_d(X;Z)$ and define the associated mesh matrix for $d$ boundaries  $M\{b_z\}$,
a $b_d \times b_d$ integral valued matrix by
$$
M\{b_s\}_{i,j} = < b_i,b_j>.
$$

In the case of $X=\Gamma$, $d=1$,  a connected graph, the study of the  mesh matrix for 1 cycles, i.e., circuits, was initiated
 by Trent \cite{Trent}, who showed
that in this case:
$$
det(M\{z_r\}) = \# \ of \ spanning / trees \ of \ the \ graph \ \Gamma
$$
Of course, the  determinant of  $M\{z_r\}$ is the product of the  eigenvalues of the real symmetric
matrix $M\{z_r\}$.

\vspace{.1in}
Here Trent's analysis of the determinant
  is extended to each dimension of any finite  CW-complex and  to
the coefficients of the complete
characteristic polynomial, $det( t \ Id - M\{z_r\})$  of the mesh matrix for $d$ cycles $M\{z_r\}$.
Also, parallel results for $M\{b_r\}$, the mesh matrix for $d$ boundaries, will also be obtained.

\vspace{.1in}
Note that under change of basis of the $d$-cycles $z_j \mapsto \Sigma_k \ A_{j,k} z_k$
with $A \in GL(z_d,Z)$, the mesh matrix $M\{z_r\}$ changes to the mesh matrix to $A \cdot M\{z_r\} \cdot A^t$,
with thus  $det( t \ Id - M\{z_r\})$  changing to $det( t \ Id - A \cdot M\{z_r\} \cdot A^t)$.
( It is only the determinant, $det(M\{z_r\})$
which is left invariant under all transformations $ M\{z_r\} \mapsto A \cdot M\{z_r\} \cdot A^t$.This last holds since $A \in GL(z_d,Z)$ implies that $det(A) = \pm 1$.)
 So it
is remarkable that one can arrive at uniform geometric formulas for
coefficients of the characteristic polynomial recording the dependence on the choice
of integral basis used.
 This same
observation also applies to $det( t \ Id - M\{b_r\})$, the characteristic matrix
for the mesh matrix for $d$ boundaries.

\vspace{.1in}
Both  mesh matrices have a useful  decomposition:
Let $A\{z_r\}, A\{b_s\}$, respectively,  denote the integral valued $n_d \times z_d,n_d \times b_d$ matrices representing the
inclusions
$$
i_z \ : \ Z_d(X;Z) \subset C_d(X;Z)\ and \ i_b  : \ B_d(X;Z) \subset C_d(X;Z)
$$
with respect to the chosen  integral  cycle basis $\{z_r\}$ of $Z_d(X;Z)$ and
 the chosen integral  boundary basis $\{b_s\}$ of $B_d(X;Z)$ and   the cellular basis of $C_d(X;Z)$.

Then the respective  square $z_r \times z_r$ mesh matrix $M\{z_r\}$ and the $b_d \times b_d$ matrix $M\{b_r\}$
are given by:
$$
M\{z_r\} = A\{z_r\}^t \cdot A\{z_r\} \ and \ M\{b_s\} = A\{b_s\}^t \cdot A\{b_s\}
$$

More intrinsically, endow $C_d(X;R)$, $Z_d(X;R)$ , $B_d(X;Z)$ with the inner products
for which the specified basis elements, cellular, d cycle, d boundary, respectively,  have length one and are orthogonal.
Denote by $<.,.>, , <.,.>^{\#_z}, <.,. >^{\#_b}$ the respective inner products.

Let  $i_z^\star :C_d(X;R) \rightarrow Z_d(X;R)$  the adjoint of $i_z$, and
$i_b^\star :C_d(X;R) \rightarrow B_d(X;R)$ the adjoint of $i_b$  be  defined by :
$$
\begin{array}{l}
<i_z(a),b> = < a, i_z^\star(b)>^{\#_z} \ for \ a \in Z_d(X;R), \ b\in C_d(X;R) \\
<i_b(a),b> = < a, i_b^\star(b)>^{\#_b} \ for \ a \in B_d(X;R), \ b\in C_d(X;R) \\
\end{array}
$$
Then the entries of $M  \{ z_r \}  $ are just
$$
M \{z_r\}_{u,v} =  <i_z(z_u), i_z(z_v)> = < z_u, i_z^\star \cdot i_z (z_v)>^{\#_z}.
$$
showing   that the mesh matrix for $d$ cycles $M\{z_r\}$
is just the matrix representing the linear mapping:
$
i_z^\star \cdot i_z : Z_d(X;R) \rightarrow Z_d(X;R)
$
and  is  here given in terms of  $A\{z_r\}$
as described above. Similarly, the mesh matrix for $d$ boundaries  $M\{b_r\}$
is just the matrix representing the linear mapping:
$
i_b^\star \cdot i_b : B_d(X;R) \rightarrow B_d(X;R)
$
and  is thus given in terms of  $A\{b_s\}$
as described above.

\vspace{.1in}

Also, one may use the cellular basis, or equivalently the inner product $<.,>$,
to define the adjoint $\partial^\star_d: C_{d-1}(X;Z) \rightarrow C_{d}(X;Z)$
 of the boundary mapping on integral $d$  chains:
$$
\partial_d :C_{d}(X;Z) \rightarrow C_{d-1}(X;Z)
$$
If $K$ represent the boundary mapping  in the cellular basis, the combinatorial
Laplacian acting on (d-1)-chains is
$$
\Delta_{d-1}= \partial_d \cdot \partial_d^\star :C_{d-1}(X;Z) \rightarrow C_{d-1}(X;Z).
$$
It is represented by the matrix $K \cdot K^t$.
 Hence, this setting
has a parallel decomposition to  that for the mesh matrices for
cycles and boundaries. This differs slightly from the standard
usage of this term in differential geometry, as explained in the
next paragraph.

Note: Logically following differential geometry, one could for $d-1>0$
follow the convention that the Laplacian on $d-1$ chains is
$\tilde{\Delta}_{d-1} = \partial_{d} \cdot \partial_d^\star + \partial_{d-1}^\star \partial_{d-1}$
where $\partial_{d-1}: C_{d-1}(X;R) \rightarrow C_{d-2}(X;R)$.
However, in view of the combinatorial Hodge decomposition, $C_{d-1}(X;R) = \mathcal{H}_{d-1} \oplus
\{ image \ \partial_{d} \} \oplus \{ image \ \partial_{d-1}^\star\}$,
an orthogonal direct sum decomposition, the mapping $\partial_{d}^\star$
maps the eigenvectors  of  $(\partial_{d}\partial_d^\star| \{image\ \partial_{d}\})$
to eigenvectors of
$(\partial^\star_d\partial_{d}| \{image\ \partial_d^\star\})$
with the same non-zero eigenvalues. Since the harmonics $\mathcal{H}_d$ are identified with
 $ \{ c \in C_d(X;R) \ | \ \partial_d (c)=0, \partial_d^\star(c) =0\}$
 and under the combinatorial   deRham isomorphism map isomorphically to the $d^{th}$ homology,
 $H_d(X;R)$, to record the non-zero eigenvalues of $\tilde{\Delta}_{d-1}$ for all $d$ efficiently
 it suffices to record the non-zero eigenvalues of what is called
 here the combinatorial Laplacian $\Delta_{d-1}= \partial_{d} \cdot \partial_d^\star$ acting on (d-1)-chains.
 This observation was used in the work of D. Ray and I. M. Singer \cite{RaySinger}
 on their  analytic version of Reidemeister-Franz torsion.

In the context of a simple graph, say $\Gamma$, $d=1$,   $C$  is just   the  familiar incidence matrix from
edges to vertices. Moreover,
$$
Kirchhoff\ matrix = Graph \ Laplacian \ \Delta_0 \ on \ 0-chains  = C \cdot C^t= \partial_1 \partial_1^\star = Deg - Adj
$$
where $Adj$ is the adjacency matrix with $Adj_{v,w} = +1$ if $\{v,w\}$ is an edge
and $=0$ otherwise, while   $Deg$ is the diagonal  matrix with $Deg_{v,v}$ equal to the
number of edges containing the vertex $v$.

\vspace{.1in}

For  any subset $V \subset S_d$  of the $d$ cells, form the d dimensional subcomplex of $X$, called $X_V$
as follows:
$$
X_V = X^{(d-1)} \cup V
$$
where $X^{(d-1)}$ is the $d-1$-skeleton of $X$, the union of all cells
of dimension at most $d-1$. By definition, $X_V$ has dimension $d$,
with $|V|$ cells of dimension $d$, and the already oriented cells of $V$
provide an integral, real  basis for $C_d(X_V;Z), C_d(X_V;R)$  respectively.

For $V \subset S_d$, let $V^c$ denote the complement of $V$ in $S_d$, \ $V^c = S_d -V$.

 \vspace{.1in}
 Now consider the boundary mapping
 $$
 \partial_d : C_d(X;Z) \rightarrow C_{d-1}(X;Z)
 $$
 with image the free group $B_{d-1}(X;Z)$ of rank $b_{d-1}$. Suppose $b_{d-1} >0$.
 By adding successively  one by one $d$-cells of $S_d$ to the $d-1$ skeleton of $X$, $X^{d-1}$,
 the rank the image under $\partial_d$ may be increased one by one until the maximal
 value $b_{d-1}$ is obtained. That is, one may exhibit subsets $V \subset S_d$
 with
 $$
 |V| = b_{d-1} \ and \ C_d(X_V,R) \stackrel{\partial_d}{\rightarrow} B_{d-1}(X_V;R) \subset B_{d-1}(X;R) \ an \ isomorphism.
 $$
In agreement with the graph case, one calls such a choice of $V \subset S_d$ with these two
 properties a \textbf{spanning forest} of dimension $d$
of $X$.  The kernel of the composite for such $V$ is just $Z_d(X_V;R)$, so alternatively,
$V \subset S_d$ is a (d-dimensional) spanning forest  precisely if
$$
|V| = b_{d-1} \ and \ Z_d(X_V;R) =0.
$$

The computation  of the determinant of the mesh matrix in dimension $d$, $det(M\{z_r\})$, will be as a sum over spanning
forests of dimension $d$ in $X$.

\vspace{.1in}
Correspondingly, given a spanning forest $V$, for each oriented d-cell, say $\sigma \in V^c$, the boundary $\partial_d \sigma$
lies in $B_{d-1}(X;R) = B_{d-1}(X_V;R)$ so is of the form $\partial_d c$ for a unique $c \in C_d(X_V;Q)$, by $V$ a spanning forest of $X$.
In this way, a unique rational d-cycle $z[\sigma]\in Z_d(X_{V \cup \{\sigma\}};Q)$ is defined such that
its restriction to $\sigma$ has coefficient $+1$. In particular, this shows that
if $\{\sigma_1,\sigma_2, \cdots, \sigma_k\}$ are disjoint oriented cells of $S_d$ in $V^c$
[with necessarily $k \le |V^c| = n_d-b_{d-1} = z_d$], the union
$$
W = V \cup \{\sigma_1,\sigma_2, \cdots, \sigma_k\},  \ 0 \le k \le z_d
$$
has the two properties
$$
|W| = b_{d-1} + k \ and\ Z_d(X_W;R) = k.
$$

One calls such a set $W \subset S_d$ a \textbf{k augmented spanning forest} since it arises
exactly by adding k d-cells to some  (d-dimensional) spanning forest $V$.
In particular, a spanning forest is just a 0 augmented spanning forest  in this notation.

Recall of the short exact sequence,
$$
0 \rightarrow Z_d(X;R) \subset C_d(X;R) \stackrel{\partial_d}{\rightarrow} B_{d-1}(X;R) \rightarrow 0.
$$
Again, in view of this  short exact sequence, an equivalent condition
for $W$ to be a k augmented spanning forest  is that $  0 \le k \le z_d$ and
$$
|W^c | = z_d -k \ and \ j_W : Z_d(X;R) \subset C_d(X;R) \stackrel{restrict}{\rightarrow } C_d(X_{W^c};R) \ is \ onto.
$$

For example,
  if $W$ is a k augmented spanning forest  of dimension $d$, then  $|W^c|  =n_d - (b_{d-1}+k) = z_d-k$.
 Since the mapping from $Z_d(X;R)$ to $C_d(X_{W^c};R)$ has kernel $
 Z_d(X;R) \cap C_d(X_W;R) = Z_d(X_W;R)$
 of dimension $k,$
   the image of $Z_d(X;R)$ in $C_d(X_{W^c};R)$   has dimension
  $ z_d - k $. Since $|W^c| = z_d -k$, this implies that $Z_d(X;R) \rightarrow C_d(X_{W^c};R)$
 is onto, proving that the first definition implies the second.
 The converse follows in a similar manner.

The computation  of the coefficients of the  characteristic polynomial of the mesh matrix
 for  $d$ cycles
 of $M\{z_r\}$, will be in terms of a sum over (d-dimensional)  k augmented spanning forests
 with suitable positive integral weights.

 \vspace{.1in}
For a finitely generated  abelian group $G$ let $[G]$ denote its torsion subgroup and $|[G]|$
the order of this subgroup.

For a finite cell complex $X$, let $t_k(X)$ denote the order of the
torsion subgroup of the $k^{th}$ homology $H_k(X;Z)$, i.e. ,
$t_k(X) = | [ H_k(X;Z)]|$.

\vspace{.3in}
By definition,  for $W$ a k augmented  spanning forest, the mapping
$$
j_W : Z_d(X;Z) \subset C_d(X;Z) \rightarrow C_d(X,X_W;Z) = C_d(X_{W^c};Z)
$$
is rationally onto. Now using the chosen cycle basis for $\{z_r\}$ for $Z_d(X;Z)$ and
the cellular basis for $C_d(X_{W^c};Z)$, the adjoint is a well defined injection
$$
j_W^\star : C_d(X_{W^c};Z) \rightarrow Z_d(X;Z) .
$$
More formally, take the cellular metric $<.,.>$ restricted to $C_d(X_{W^c};Z)$ and
the cycle metric  $<.,.>^{\#_z}$ on $Z_d(X;Z)$ and form the adjoint. Hence,
if  $j_W$ is represented in these basis by the matrix, say $U_W$, an integer
valued $(z_d-k) \times z_d$ matrix, then the adjoint is represented by
$U_W^\star$, the transpose.

The determinant $det(j_W \cdot j_W^\star)$
will be called the \textbf{covolume squared} of this inclusion $j_W^\star$ of the lattice.
It may be seen to  be equal to the square of the unit volume cell of the sublattice
inside the real vector space $Z_d(X;R)$ under the metric $<.,.>^{\%_z}$.
In   this case, as $U_W$ is an integral  valued matrix, this covolume squared is
a non-negative integer.

Similarly, one has the injection of $d$ cycles :
$$
m_W : Z_d(X_W;Z) \subset Z_d(X;Z)
$$
By $W$ a k augmented  spanning forest the  free abelian group  $Z_d(X_W;Z)$ has rank $k$. Pick
an integral  basis, say $\{\hat{z}_u\}_{u=1}^k$ of $Z_d(X_W;Z)$. For this
$d$ cycle basis and the already chosen standard basis on $\{z_r\}$ for $Z_d(X;Z)$,
let $B[W]$ be  the $z_d \times k$ matrix representing this
inclusion with respect to the  integral basis already chosen.
 As seen the determinant
$$
det( B[W]^t B[W] )  \ is \ a \ positive \ integer
$$
and is independent of the choice of integral basis for $Z_d(X_W;Z)$.
It explicitly records the dependence on the fixed basis of $Z_d(X;Z)$.

This determinant is   employed in the first main theorem.

\begin{thm} \label{thm1}
The constant coefficient of $(-1)^{z_d} \ det(t\ Id - M\{z_r\}) $ is
$$
 \begin{array}{l}
 det( M\{z_r\})  = product \ of \ eigenvalues \ of \ M\{z_r\} = \sigma_{z_d}(M\{z_r\})  \\
    = \Sigma_{ V, \ a \ spanning \ forest \ of \ dimension \ d}   \ \  (\frac{t_{d-1}(X_V)}{t_{d-1}(X)} )^2 \ . \\
 \end{array}
 $$
Also, for $V$, any  spanning forest, i.e, a 0 augmented spanning forest,
$$
 \ covolume \ squared \ of  (\ j_V^\star : C_d(X_{V^c};Z) \hookrightarrow Z_d(X;Z) \ ) =
 det( j_V \cdot j_V^\star) =  (\frac{t_{d-1}(X_V)}{t_{d-1}(X)})^2 \ .
 $$

If $ 0  <  k < z_d$, then
$$
 \begin{array}{l}
 \  coefficient \ of \ (-1)^{z_d-k} \  t^{k} \ in \ det(t \ Id - M\{z_r\}) \\
= (z_d-k)^{th}  \ elementary \ symmetric \ function \ in \ eigenvalues \ of \ M\{z_r\} = \sigma_{z_d - k}( M\{z_r\}) \\
   = \Sigma_{ W, \ a \ k \ augmented \  spanning \ forest }  \ (\frac{t_{d-1}(X_W)}{t_{d-1}(X)})^2 \
   det( B[W]^t B[W]) \ .
 \end{array}
 $$
 Here the integral matrix $B[W]$ records, as above,  the injection $m_W : Z_d(X_W;Z) \subset  Z_d(X;Z)$.

Moreover, for $W$, a k augmented spanning forest, with $0 \le k< z_d$
$$
\begin{array}{l}
 \ covolume \ squared \ of  (\ j_W^\star : C_d(X_{W^c};Z) \hookrightarrow Z_d(X;Z) \ ) = det( j_W \cdot j_W^\star) \\
 =  (\frac{t_{d-1}(X_W)}{t_{d-1}(X)})^2 \cdot det( B[W]^t B[W]) \ .
  \end{array}
 $$

Additionally, the quotient $ \frac{t_{d-1}(X_W)}{t_{d-1}(X)}$ is a positive  integer for $W$
any k-augmented spanning forest;  in particular, for $k=0$.
\end{thm}

The first part extends Trent's theorem to cell complexes; the second extends it
 to evaluate the complete characteristic
polynomial of the mesh matrix $M\{z_r\}$. Note that
two separate expressions are given for the  integer weights,
one involving $X_W$ and the other $X_{W^c}$. The explicit dependence
on the choice of the integer basis $\{z_r\}$ made for $Z_d(X;Z)$ is recorded directly through
$det(B[W]^t B[W])$.

 In \S \ref{sectiongeombases}, for a geometrically
defined basis of $Z_d(X;Q)$, given by a choice of d dimensional spanning forest,$V_0$,
the coefficients of the complete characteristic polynomial
of the cycle mesh matrix associated to that geometrical choice
are given combinatorially. Similarly, in \S \ref{sectiongeombases},  for
 a (d+1) dimensional spanning forest $V_1$, there are geometrically defined bases for $B_d(X;Q)$.
Again in that section, the complete characteristic polynomial of the
associated geometrically defined cycle and boundary mesh matrices are  evaluated combinatorially.

\section{\textbf{Proof of the main theorem for mesh matrix for $d$ cycles: Theorem \ref{thm1}}}  \label{mesh1}
\vspace{.1in}

In analyzing $M\{z_r\} = A\{z_r\}^\star \cdot A\{z_r\}$ and its characteristic polynomial
 the following  elementary standard lemmas are helpful.

 Let $A$ be a $m \times n$ matrix. For each pair $I,J$ of subsets $I \subset \{1,2, \cdots, m\}$, $J \subset \{1,2, \cdots, n\}$
 let $A_{I,J}$ be the $|I| \times |J|$  minor of $A$ with row indices from the subset $I$ and column indices from the subset $J$.

\begin{lemma} \label{lemmaalg1}
Let $A$ be square $n \times n$ matrix, then $\sigma_k(A)$ the $k^{th}$, elementary
symmetric function of the eigenvalues, say $\lambda_i, i =1 ,\cdots, n$, counted with
multiplicities, is given by the sum of determinants
$$
\sigma_k(A) = \Sigma_{I \subset \{1,\cdots, n\} \ with \ |I| =k } \   \ det(A_{I,I})
$$
for $ 0 < k \le n$.
\end{lemma}

\begin{lemma} [via Cauchy-Binet theorem.] \label{lemmaalg2}

Let $A,B$ be $m \times n , n \times m$ matrices respectively, then  the non-vanishing
eigenvalues of $A \cdot B$, a $m \times m$ matrix, and $B \cdot A$, a \ $n \times n $ matrix
counted with multiplicities are identical and the $k^{th}$ elementary symmetric function of the eigenvalues of $AB$
or of $BA$ equals:
$$
\begin{array}{l}
\sigma_k(AB)
= \Sigma_{I \subset \{1, \cdots, m\} \ with |I| =k } \ det((A B)_{I,I}) \\
= \Sigma_{I \subset \{1, \cdots, m\} \ and \  J \subset \{1,\cdots, n \} \ with \ |I|=k, |J|=k} \ det(A_{I,J}) det(B_{J,I})
\ \ Cauchy-Binet \\
= \Sigma_{I \subset \{1, \cdots, m\} \ and \  J \subset \{1,\cdots, n \} \ with \ |I|=k, |J|=k} \ det(B_{J,I}) det(A_{I,J}) \\
= \Sigma_{J \subset \{1, \cdots, n\} \ with |J| =k } \ det((A B)_{J})  \ \ Cauchy-Binet\\
 = \sigma_k(AB) \\
 \end{array}
$$

\end{lemma}

To prove the first lemma is elementary. The desired coefficient may be evaluated
by differentiating $det(t \ Id - A)$ a total of $(n-k)$ times and setting $t=0$.
Due to the special form of the square matrix $(t \ Id -A)$, with $t's$ only down the diagonal, this
gives a sum of determinants of minors along the diagonal blocks of size $(k \times k) = |I| \times |I| $ indexed by the sets $(I,I)$.

The second lemma follow from the first by utilizing, as in the statement of the theorem,
the Cauchy-Binet theorem.

\vspace{.1in}

These lemmas can be applied to the mesh matrix on $d$ cycles $M\{z_r\} = A\{z_r\}^t \cdot A\{z_r\}$ for $d$ cycles
where  $A\{z_r\}$, a $n_d \times z_d$ integer valued matrix, represents the inclusion:
$$
i_X : Z_d(X;Z) \subset C_d(X;Z).
$$
The rows of $A\{z_r\}$ are indexed by d-cells in $S_d$ and the columns by indices in $\{1,\cdots, z_d\}$
specifying  elements of the cycle basis $\{z_r\}$ of $Z_d(X;Z)$.

By the second  lemma the $(z_d-k)^{th}$ elementary symmetric function in the eigenvalues
of $M\{z_r\} = A\{z_r\}^t \cdot A\{z_r\}$ equals the $(z_d-k)^{th}$ elementary symmetric function in the eigenvalues
of  $A\{z_r\} \cdot A\{z_r\}^t$. By the first  lemma this equals the sum of determinants,
$$
\Sigma_{I \subset S_d \ with \ |I| = z_d -k} \ det( \ (A\{z_r\} \cdot A\{z_r\}^t)_{(I,I)}\ ),
$$
where the sum, as indicated, is over the $d$-cells $S_d$ of $X$ which label the
rows of $A\{z_r\}$.

Here there is the equalities of matrices,
$$
\begin{array}{l}
(A\{z_r\} \cdot A\{z_r\}^t)_{(I,I)}
= (A\{z_r\})_{(I,\{1,\cdots, z_d\}) }   \cdot (A\{z_r\}^t)_{(\{1, \cdots,z_d\} ,I)} \\
= (A\{z_r\})_{(I,\{1,\cdots, z_d\}) } \cdot  ( \ (A\{z_r\})_{(I,\{1,\cdots, z_d\}) } \ )^t\ , \\
\end{array}
$$
where $(A\{z_r\})_{(I,\{1,\cdots, z_d\})}$ represents
$$
r_{I} : Z_d(X;Z) \subset C_d(X;Z) \stackrel{restrict}{\rightarrow}  C_d(X_I;Z)
$$
and  its transpose the adjoint for the cycle basis for $Z_d(X;Z)$ and the cellular basis
for $C_d(X_I;Z)$. In toto, this proves:
$$
\begin{array}{l}
The\ (z_d-k)^{th} \ elementary \ symmetric \  function \ in \  the  \ eigenvalues\
of \ M\{z_r\}  \\
=\Sigma_{I \subset S_d \ with \ |I| = z_d -k} \ det( r_{I} \cdot r_{I}^\star) \ .
\end{array}
$$

\vspace{.1in}
Now clearly, the determinant $det(r_{I} \cdot r_I^\star)$ equals zero unless \newline
$r_I : Z_d(X;Z) \subset C_d(X;Z) \stackrel{restrict}{\rightarrow}  C_d(X_I;Z)$ is onto.

Here $|I| = z_d-k$, so this is precisely the condition
that $I=W^c$ with $W$ a k augmented spanning forest. In this case,
$$
r_{I} = j_W : Z_d(X;Z) \subset C_d(X;Z) \stackrel{restrict}{\rightarrow}  C_d(X_{W^c};Z) \ .
$$

That is, this proves:
$$
\begin{array}{l}
^he\ (z_d-k)^{th} \ elementary \ symmetric \  function \ in \  the  \ eigenvalues\
of \ M\{z_r\}  \\
=\Sigma_{W \subset S_d \ a \ k \ augmented \ spanning \ forest  } \ det( j_{W} \cdot j_{W}^\star) \ .
\end{array}
$$
Here the adjoint of $j_W$ called here $j_W^\star$ is formed by declaring that the
fixed chosen basis of $Z_d(X;Z)$ is orthonormal and that the d-dimensional cellular basis of
$ C_d(X_{W^c};Z)$ is orthonormal.

\vspace{.1in}
To complete the proof of the theorem for the $k=0$ case, for a spanning forest, say $V$,
one needs to show that
the covolume squared is expressed as
$$
det(j_V \cdot j_V^\star) = (\frac{t_{d-1}(X_V)}{t_{d-1}(X)})^2 \ .
$$

Now consider
 $$
 j_V : Z_d(X;Z) \rightarrow C_d(X_{V^c};Z).
 $$
It is   a  mapping
 of free abelian groups and is rationally
an isomorphism by the definition of spanning forest,  so $det(j_V \cdot j_V^\star)$
is just the square of the order of the torsion abelian group $C_d(X_{V^c};Z)/j_V(Z_d(X;Z))$.

Now quite generally, let $W$ be a k-augmented spanning tree and consider the restriction mapping
$$
j_W: Z_d(X;Z) \rightarrow C_d(X_{W^c};Z)
$$
which is rationally onto by $W$ an k-augmented spanning forest, it
is claimed that  the torsion group
$
C_d(X_{W^c};Z)/j_V(Z_d(X;Z))$
has order
$$
| \ C_d(X_{W^c};Z)/j_W(Z_d(X;Z)) \ | = \frac{t_{d-1}(X_W)}{t_{d-1}(X)}\ ,
$$
in particular is an integer as asserted in theorem \ref{thm1}.

Hence, in the case $k=0$, this implies the desired formula above.

 \vspace{.3in}
To see this last write $j_W$ as the composite
$$
j_W :Z_d(X;Z) \subset  C_d(X;Z) \rightarrow C_d(X_{W^c};Z) =  C_d(X,X_W;Z) = Z_d(X, X_W;Z)
$$
since $X$ and $X_W$ have the same (d-1) skeleton. Also note that $B_d(X;Z)$ maps
onto $B_d(X,X_W;Z)$, so the desired quotient is isomorphic to the quotient
$$
H_d(X,X_W;Z)/\{ image \ H_d(X;Z) \} \ .
$$

But recall  from Lyons \cite{Lyons} the utility of introducing  the
 long exact sequence for the pair $(X, X_W)$ which reads
$$
H_d(X;Z) \stackrel{\alpha}{\rightarrow} H_d(X, X_W;Z) \stackrel{\beta}{\rightarrow}
 H_{d-1}(X_W;Z) \stackrel{\gamma}{\rightarrow} H_{d-1}(X;Z) \rightarrow H_{d-1}(X,X_W;Z)=0
$$
with $\alpha$ is   rationally onto by $W$ a k augmented spanning forest.

Hence, letting $K$ denote the image of $\beta$ in $H_{d-1}(X_W;Z)$, by exactness,
the torsion group $C_d(X,X_W;Z)/j_V(Z_d(X;Z)) \cong H_d(X,X_W;Z)/\{ image \ H_d(X;Z)\}$ is
isomorphic to \newline  $ image \ of \ \beta  = K$.

 Now consider the associated mappings of torsion subgroups:
 $$
 image \ \beta = K =   [K]  \subset  [H_{d-1}(X_W;Z)] \stackrel{\gamma'}{\rightarrow} [H_{d-1}(X;Z)] \rightarrow 0
  $$
By definition $\gamma' \cdot \beta =0$. Now let  $x$  be a torsion class
 in $H_{d-1}(X;Z)$, say $Nx =0, $ with $N$ a positive integer.
By exactness of the pair, there is $y \in H_{d-1}(X_W;Z)$ with $\gamma(y) = x$. But by
$\gamma(Ny) =0$, $Ny$ is in the image of $\beta$ which has only torsion image.
In particular,  $Ny $ is  torsion and so also is $y$. Hence, $\gamma'$ is onto. Also, this shows the sequence is exact
at $[H_{d-1}(X_W;Z)]$. In toto, this proves the desired equality
$$
 |K|  = | [H_{d-1}(X_W;Z)]| /|[ H_{d-1}(X;Z)]|  =  \frac{t_{d-1}(X_W)}{t_{d-1}(X)} = a\ positive  \ integer \ .
$$

The remaining parts of theorem \ref{thm1} will be addressed in the next section \S \ref{sectionalter}.

\section{\textbf{Alternative version of the squared covolume $det(j_W \cdot j_W^\star)$, cases $ k > 0$}} \label{sectionalter}
\vspace{.1in}

To complete the proof of theorem \ref{thm1}, one must prove
$$
squared \ covolume \ det(j_W \cdot j_W^\star) =
(\frac{t_{d-1}(X_W)}{t_{d-1}(X)})^2 \ det(B[W]^t B[W])
$$
for $W$ a k-augmented spanning forest, $0 < k < z_d$.

To prove this result the following
algebraic lemma will play a crucial part. It is proved at the end of this section.

\begin{lemma} \label{lemma55}

Let $1 \le u < n$ and $A$ be an invertible  n by n matrix. Let $R$ be the
$u \times n$ matrix obtained from $A$ by deleting the last $n-u$ rows of $A$.
Let $S$ be the matrix obtained from $(A^{-1})^t$, the transpose inverse of $A$
by deleting the  first $u$ rows, then
$$
  det(R \cdot R^t) = det(A)^2 \cdot det( S \cdot S^t) \ .
$$
\end{lemma}

To be  more explicit, as before let $\{z_r\}$ be the fixed chosen basis for $Z_d(X;Z)$,
and
$$
i_X : Z_d(X;Z) \subset C_d(X;Z)
$$
be the inclusion.
Use the \textbf{cellular basis},  $\{1 \cdot \sigma_j \ | \ j = 1 , \cdots, n_d\}$ of oriented d cells
to define the $n_d \times z_d$ matrix $A$ by
$$
i_X(z_r) = \Sigma_{\sigma_t \in S_d} \ A_{t,r} \ \sigma_t \ .
$$
\{Here, as is conventional, $i_X(\Sigma_r L_r z_r) = \Sigma_k( \Sigma_r A_{k,r} L_r) \sigma_k$.\}

Now let $W$ be a k-augmented spanning forest, say $W = V \cup I$,  with $V$ a choice of spanning
tree and $I \subset S_d -V$ of size $k$.

By $V$ a spanning tree,  there is an isomorphism
$$
j_V : Z_d (X;Q) \rightarrow C_d(X_{V^c};Q)
$$
of rational vector spaces.

 Indeed,  for each cell $\sigma_j \in V^c$
there is a unique rational d-cycle, say $z( \sigma_j) $ for $X$
which maps under $j_V$ to precisely the d-chain $1 \cdot \sigma_j$. That is,
$$
z(\sigma_j) \in Z_d( X_{V\cup \sigma_j};Q) \ restricts \ to \ 1 \cdot \sigma_j
$$
for each $\sigma_j \in V^c$. This geometrically defined cycle $z(\sigma_j), j=1 ,\cdots, z_d$,
with $\sigma_j$ ranging over $V^c$,
is uniquely determined by these requirements.

In  this manner, having chosen the spanning forest $V$,  one obtains a canonical rational
geometric basis $\{z(\sigma_j)\ | \ \sigma_j \in V^c= S_d -V\}$ for $Z_d(X;Q)$ which restricts under $j_V$ to precisely the
cellular orthonormal basis of $C_d(X_{V^c};Q)$. Very nicely the associated natural geometric
basis for $Z_d(X_W;Q)$ with the same spanning forest $V \subset W$ is precisely
$\{ z(\sigma_j) \ | \ \sigma_j \in W - V \}$ and similarly  $Z_d(X_{W^c};Q)$ has rational basis
$\{ z(\sigma_j) \ | \ \sigma_j \in W^c\}$. In this manner, this geometric basis
respects the direct sum decomposition of $Z_d(X;Q)$: [$W = V \cup I$,$ S_d - W = V^c \cup I$ $|I| =k$.]
$$
Z_d(X;Q) = Z_d(X_W;Q) \oplus Z_d(X_{V \cup W^c)};Q) = Z_d(X_{V \cup I};Q) \oplus Z_d(X_{V \cup W^c)};Q)
$$
In particular the methods of lemma \ref{lemma55} apply for this choice of geometrical basis,
those for $\sigma_j \in I$ and those with $\sigma_j \in W^c$.

As has been demonstrated above in terms of the restriction mapping $j_X : Z_d(X;Z) \rightarrow C_d(X_{V^c};Z)$
which is rationally an isomorphism with  respect to
an integral basis of $Z_d(X;Z)$ to the cellular one has  determinant $\pm \frac{t_{d-1}(X_V)}{t_{d-a}(X)}$.
Alternatively expressed, the change of basis from the an integral basis of $Z_d(X;Z)$
to the geometrical basis given by the spanning forest $V$
has determinant equal to $\pm \frac{t_{d-1}(X_V)}{t_{d-1}(X)}$.
As observed
above, the $k=0$ case of theorem \ref{thm1}
$
 det(A) = \pm (\frac{t_{d-1}(X_V)}{t_{d-1}(X)}).
 $

In a like manner, since  $V$ as a spanning forest for $X_W$, the
 change of basis from the an integral basis of $Z_d(X_W;Z)$
to the geometrical basis given by the spanning forest $V$ in $X_W$
for $Z_d(X; Q)$
has determinant equal to $\pm \frac{t_{d-1}(X_V)}{t_{d-1}(X_W)}$.

\vspace{.3in}

Under this identification of $Z_d(X;Q)$ with $C_d(X_{V^c};Q)$, the mapping
$j_V$ becomes identified with the identity mapping
$$
Id : Z_d(X;Q) \cong  Z_d(X;Q)
$$
with the  the matrix $A$ becoming just the translation between the two bases
$\{z_r\}$ and the geometric basis $z(\sigma_j)$ as in
$$
z_r = \Sigma_{ 1 \le j \le z_d} \ \ A_{j,r} \ z(\sigma_j) \ .
$$
[As is seen by applying the isomorphism $j_V$ to both sides.]

 Now arrange that the first $z_d-k$ cells of $V^c$ are those of $W^c=S_d - W = S_d - (V \cup I)= V^c -  I $
 a total of $z_d-k$ in all.  Let $R$ be  the $(z_d-k) \times z_d$
 matrix obtained from $A$ by deleting the last $k$ rows. This matrix represents
 exactly the mapping
 $$
 j_W : Z_d(X;Q) \stackrel{restrict}{\rightarrow} C_d(X_{W^c};Q) \cong Z_d(X_{V \cup W^c};Q)
 $$
 sending $z_r \mapsto  \Sigma_{1 \le j \le z_d-k}\ A_{j,r} \ z(\sigma_j)$
 [ as seen by applying $j_W$]
 and $det(j_W \cdot j_W^\star)$ is represented by the matrix $R R^T$.

 On the other hand, inverting the matrix $A$, gives
 $$
 z(\sigma_t) = \Sigma_r \ (A^{-1})_{r,t} \ z_r \ .
 $$
 That is, the identify mapping
 $$
 Id : Z_d(X;Q) \cong Z_d(X;Q)
 $$
 which for the geometric basis on the domain, $z(\sigma_j)$ with $ \sigma_j \in V^c$
and the fixed chosen basis $\{z_r\}$ on the range is represented by
the matrix $A^{-1}$.

Hence, the inclusion of d-cycle spaces over the rationals
$$
i_W : Z_d(X_W;Q) \subset Z_d(X;Q)
$$
is represented for the domain by the geometrical basis $z(\sigma_j), 1 \le j\le k$ $\sigma_j \in I$,
and range basis $\{z_r\}$ by the matrix, say $S^t$, and is obtained from $A^{-1}$ by
 deleting the first $k$ columns. Equivalently, $i_W^\star$ is represented by $S$
 which is obtained form $(A^{-1})^t$
by deleting the first $k$ rows.

Since $B[W]$ by definition is the matrix for the inclusion $Z_d(X_W;Z) \subset Z_d(X;Z)$
with respect to the an integral basis for $Z_d(X_W;Z)$ and the chosen
basis $\{z_r\}$ for $Z_d(X;Z)$, this observation shows that
$$
 det(B[W]^t \cdot B[W]) = ( \frac{t_{d-1}(X_V)}{t_{d-a}(X_W)})^2 \cdot det( i_W^\star \cdot i_W)
$$
and $i_W \cdot i_W^\star$ is represented by
$$
\begin{array}{l}
 ( \  (A^{-1})^t_{r,t} |_{1 \le  t \le z_d-k} \ )
 \cdot (A^{-1})_{r',t'}|_{1 \le  t' \le z_d-k}  \\
 = ( \  ((A^{-1})^{t})_{t,r} |_{1 \le  t \le z_d-k} \ )
 \cdot ( ((A^{-1})^{t})^t_{t',r'}|_{1 \le  t' \le z_d-k} )^t \\
 = S \cdot S^t \ .
 \end{array}
$$

By  lemma \ref{lemma55}:
$$
det(R \cdot R^t) = det(A)^2 \cdot det(S \cdot S^t) \ .
$$

  Hence, one  sees that,
in  the notation of theorem \ref{thm1}:
$$
\begin{array}{l}
det( \ j_{W}^\star \cdot j_{W}\ ) = det( R \cdot R^t)
= (det(A))^{2} \cdot det(S \cdot S^t) \\
= (\frac{t_{d-1}(X_V)}{t_{d-1}(X)})^2 \cdot ( \frac{t_{d-1}(X_V)}{t_{d-1}(X_W)})^{-2} \cdot det(B[W]^t B[W])
= (\frac{t_{d-1}(X_W)}{t_{d-1}(X)})^{2} \ \cdot det(B[W]^t B[W])
\end{array}
$$

As noted above, $B[W]$ is an integral valued matrix and $\frac{t_{d-1}(X_W)}{t_{d-1}(X)}$
is an integer when $W$ is a k augmented spanning forest. So this completes the proof of
theorem \ref{thm1}.

\vspace{.3in}
\textbf{Proof of lemma \ref{lemma55}}

One first notes that by the Cauchy-Binet theorem
the determinant of $R \cdot R^t$ equals  the sum of
squares of the $r \times r$ minors of $A$
 which occur in
the first $r$ rows. Similarly, the determinant of $S \cdot S^t$ equals  the sum of
squares of the $s \times s$ minors of $(A^{-1})^t$  which occur in
the last  $s$ rows. Hence, the above lemma will follow from the
following lemma:

\begin{lemma}
Let  $A$ be a square $(r+s) \times (r+s)$ non-singular matrix, i.e, $det(A) \neq 0$.
Let $E$ be a $r \times r$ minor of $A$ and $F$  be the complementary $s \times s$ minor
of $(A^t)^{-1}$, then
$$
det(E) = \pm det(A) \cdot det(F) \ .
$$
\end{lemma}

Now replacing $A$ by pre and post multiplication by permutation matrices,
with care about signs, this last lemma is reduced to:

\begin{lemma}
Let  $A$ be a square $(r+s) \times (r+s)$ non-singular matrix, i.e, $det(A) \neq 0$.
Let $E$ be the top left $ r\times r$ block of $A$ and $F$ be the bottom $s \times s$
block of $(A^t)^{-1}$. Then
$$
det(E) =  det(A) \cdot det(F) \ .
$$
\end{lemma}

By continuity, it suffices to consider the case with $det(E) \neq 0, det(F) \neq 0$.
For these cases this lemma  follows from the standard theorem.

\begin{thm}
Let $A$ be a square $(r+s) \times (r+s)$  non-singular matrix, of the block partitioned form:
$$
A = \left( \begin{array}{cc}    a & b  \\ c & d \end{array}   \right)
$$
where $a,b,c,c$ are the $r \times r, r \times s, s \times r , s \times s$
matrices respectively, defined by $a = ( \ A_{i,j} |_{1 \le i \le r, \ 1 \le j \le r} \ ),
b = ( \ A_{i,j}|_{1 \le i \le r, \ r+1 \le j \le r+s} \ ), c = ( \ A_{i,j}|_{r+1 \le i \le r+s,\  1 \le j \le r} \ ),
 d = ( \ A_{i,j}|_{r+1 \le i \le r+s,\  r+ 1 \le j \le r+s} \ )$ respectively.

\item If $det(a) \neq 0 $, then
$
det(A) = det(a) \cdot det( d - c \ a^{-1} \ b).
$
In particular, the inverse $ (d - c \ a^{-1} \ b)^{-1}$ exists in    the context $det(A) \neq 0, det(a) \neq 0$.

 If $det(d) \neq 0$, then
$
det(A) = det(d) \cdot det( a - b \ d^{-1}  \ c)
$
In particular, the inverse $ (a - b \ d^{-1} \ c)^{-1}$ exists in    the context $det(A) \neq 0, det(d) \neq 0$.

If $det(A) \neq 0, det(a) \neq 0, det(d) \neq 0$, then
the inverse of $A$, $A^{-1}$, has the block decomposition :
$$
\begin{array}{l}
A^{-1} = \left( \begin{array}{cc} ( a - b \ d^{-1}  c)^{-1} \  & \ -(a-b \ d^{-1} \ c)^{-1} \ b d^{-1}\\
                               - (d- c\ a^{-1} \ b)^{-1} \ c a^{-1}  \ & \ ( d - c \ a^{-1}  \ b)^{-1}
                               \end{array} \right) \\
                               \hspace{.5in}
                               = \left( \begin{array}{cc} ( a - b \ d^{-1}  c)^{-1} \  & \ - a^{-1}b \ (d - c \ a^{-1} \ b)^{-1} \\
                               - d^{-1} c \ (a- b\ d^{-1} \ c)^{-1}  \ & \ ( d - c \ a^{-1}  \ b)^{-1}
                               \end{array} \right)
\end{array}
                               $$

\end{thm}

The proof of this is elementary.
By $det(a) \neq 0$, one may reduce $A$ to
 Schure complementary form by multiplying by a matrix of determinant $1$; as in:

$$
 \left( \begin{array}{cc}    a \ \  & \ \  0  \\ c \ \  & \  \ d - c \ a^{-1} \ b  \end{array}   \right)
 = \left( \begin{array}{cc}    a \ \  & \ \  b  \\ c \ \  & \  d  \end{array}   \right)
 \cdot \left( \begin{array}{cc}    Id \ \  & \ \  -a^{-1}\ b \\ 0 \ \  & \  Id  \end{array}   \right)
$$
Hence, $det(A) = det(a) \cdot det(d - c \ a^{-1} \ c \  b)$ as claimed.
Similarly for the case $det(d) \neq 0$.

The formulas for $A^{-1}$ follow by direct calculation of $A^{-1} \cdot A$.

\section{\textbf{The mesh matrix for $d$ boundaries; evaluation of its  characteristic polynomial}} \label{sectionbd}
\vspace{.1in}

There is a parallel development for the mesh matrix on $d$ boundaries.

One fixes a choice of basis for the integral d-dimensional  boundaries $B_d(X;Z)$
say $\{b_r\}$ and declares that they form an orthonormal basis for the boundaries.

\vspace{.2in}

 Let $V \subset S_d$ be called a \textbf{spanning coforest}
 if  the restriction mapping
 $$
 k_V : B_d(X;Q) \subset C_d(X;Q) \stackrel{restrict}{\rightarrow} C_d(X_V;Q)
 $$
 is an isomorphism. In particular, $|V| = b_d= rank \ B_d(X;Q)$. This is called by Lyons \cite{Lyons} a cobase.

 By incrementally making appropriate choices of $d$-cells in $S_d$, the image of the restriction
 can be made larger and larger until it achieves the maximal value $b_d$. Hence, such
 spanning coforests always  exist.

 Now let $W$ be the result of deleting $k$ cells [ with $ 0 \le k < b_d$],
 from a spanning coforest, say $V$. The result is called a \textbf{k reduced spanning coforest}.
 It is characterized by $|W| = b_d -k$ and
 $$
 k_W : B_d(X;Q) \subset C_d(X;Q) \stackrel{restrict}{\rightarrow} C_d(X_W;Q)
 $$
 is surjective. Naturally, a
  spanning coforest is just a 0 reduced spanning coforest  in this terminology.

 This $k_W$  is obtained by tensoring with the rationals the composite mapping
 $$
 \begin{array}{l}
 k_W : B_d(X;Z) \subset C_d(X;Z) \rightarrow C_d(X,X_{W^c};Z) =C_d(X_{W};Z),
 \end{array}
 $$
 which is rationally a surjection by $W$ a k reduced spanning coforest.
 Its adjoint with respect to the orthonormal basis $\{b_r\}$ of $B_r(X;Q)$
 and the cellular basis [declared to be orthonormal] of $C_d(X_{W^c};Z)$  is
 the inclusion :
 $$
 k_W^\star : C_d(X_{W};Z) \hookrightarrow B_d(X;Z).
 $$
 Again,  one  calls the determinant $det( k_W \cdot k_W^\star)$
\textbf{the covolume squared} of this inclusion $k_W^\star$ of lattices}.

\vspace{.3in}
Analogously to the mesh cycle cases, let $\{b_s\}$ be a fixed choice of
integral basis for $B_d(X;Z)$ and introduce the boundary mesh matrix
$$
M\{b_s\} = \{  <b_i,b_j> \}
$$
where one includes $B_d(X;Z) \subset C_d(X;Z)$ and used the cellular
basis to define the inner product on d chains. Equivalently,
 this boundary mesh matrix is $k_X \cdot k_X^\star$.
 with the adjoint computed with respect to the respective basis $\{b_s\}$
 on $B_d(X;Z)$ and the cellular basis for the d chains made orthonormal
 in the respective inner products.

As seen before, the application of lemma  \ref{lemmaalg2} evaluates
the $b_d-k$ elementary symmetric function of the eigenvalues of the
boundary mesh matrix $M(\{b_s\})$ as a sum :
$$
\sigma_{b_d-k}(M(\{b_s\})) = \Sigma_{W \ with \ |W| = b_d -k} \ det( k_W \cdot k_W^\star)
$$
Again the determinant is non-zero only for the case that $k_W$ is onto, that is
the case that $W$ is a k reduced spanning coforest. Hence,
$$
\sigma_{b_d-k}(M(\{b_s\})) = \Sigma_{W \ with \ W \ a \ k \ reduced \ spanning \ coforest} \ det( k_W \cdot k_W^\star) \ .
$$

\vspace{.3in}

For  a k-reduced spanning coforest  $W$ chose a spanning coforest
$V$ with $W \subset V$. That is, $W$ is obtained from $V$ by deleting $k$ oriented cells, say $I$,
so $V = W \cup I$ and $|I|=k$.  Having made such a choice of $V$, consider the
inclusion mapping :
$$
x_W : B_d(X_{W^c};Z)  \subset B_d(X;Z)
$$
For a choice of integral basis for $B_d(X_{W^c};Z) =  B_d(X_{V^c \cup I};Z))$
and the chosen basis $\{b_s\}$ for $B_d(X;Z)$ let $B'[W]$ denote
the associated integral matrix. Hence, the $b_d \times b_d$ matrix
$B'[W]^tB'[W]$ represents $x_W^\star \cdot x_W^\star$ and $det(B'[W]^tB'[W])$
is independent of the choice of integral basis of $ B_d(X_{W^c}Z)= B_d(X_{V^c \cup I};Z))$.

By $W$ a k reduced spanning coforest equal to $W=V-I, |I| = k$ for a spanning coforest $V$, the restriction mapping
$$
y_W: B_d(X_{W^c};Z) = B_d(X_{V^c\cup I};Z)  \mapsto C_d(X_I;Z)
$$
will be a rational isomorphism.

\vspace{.1in}
The parallel result to theorem \ref{thm1}  for the mesh matrix for boundaries
 $M\{b_s\}$
is:

\begin{thm} \label{thm2}

For $V$ a d dimensional spanning coforest of $X$, the group $H_d(X;X_{V^c};Z)$
is a torsion subgroup. Let its order be denoted by $v(V,X)$.

\vspace{.3in}
The constant coefficient of $(-1)^{b_d} \ det(t\ Id - M\{b_r\}) $ is
$$
 \begin{array}{l}
 det( M\{b_s\})  = product \ of \ eigenvalues \ of \ M\{b_s\} = \sigma_{b_d}(M\{b_s\})  \\
    = \Sigma_{ V, \ a \ spanning \ coforest \ of \ dimension \ d}   \ \  (v(V,X) )^2 \ . \\
 \end{array}
 $$
Also, for $V$, any  spanning coforest, i.e, a 0 reduced  spanning coforest,
$$
 \ covolume \ squared \ of  (\ k_V^\star : C_d(X_V;Z) \hookrightarrow B_d(X;Z) \ ) = det(k_V \cdot k_V^\star) = (v(V,X))^2 \ .
 $$

\vspace{.1in}
Quite generally, let  $W$ be a d-dimensional k reduced spanning coforest in $X$.
For such $W$,  the quotient $C_d(X_{W^c};Z)/x_W(Z_d(X;Z))$
is a torsion group isomorphic to $H_d(X;X_{W^c};Z)$. Let the order of this finite torsion
group be denoted by $v(V,X)$.

Moreover for $W$  contained in a  spanning coforest $V$
 with $W = V-I$ for $|I| =k$, then  the quotient $ C_d(X_I;Z)/y_W(B_d(X_{W^c});Z))$ is a torsion group and is
 isomorphic to $H_d(X_{W^c}, X_{V^c};Z)$. Let the order of this finite torsion group be
 denoted by $u(W,V)$.

It is claimed that the quotient
$$
f(W) = \frac{u(W,V)}{v(V,X)}
$$
is independent of the choice of $V$ and moreover
$
f(W) \ v(W,X)
$
is an integer.

If $ 0 \le k < b_d$, then
$$
 \begin{array}{l}
 \  coefficient \ of \ (-1)^{b_d-k} \  t^{k} \ in \ det(t \ Id - M\{b_r\}) \\
= (b_d-k)^{th}  \ elementary \ symmetric \ function \ in \ eigenvalues \ of \ M\{b_r\} = \sigma_{b_d - k}( M\{b_r\}) \\
   = \Sigma_{ W, \ a \ k \ reduced  \  spanning \ coforest }  \ (\frac{u(W,V)}{v(V,X)})^2 \
   det( (B'[W])^t B'[W]) \ .
 \end{array}
 $$

Moreover, for $W$, a k augmented spanning forest, with $0 \ge k< z_d$ :
$$
\begin{array}{l}
 \ covolume \ squared \ of  (\ k_W^\star : C_d(X_{W};Z) \hookrightarrow B_d(X;Z) \ ) \\
 =  (\frac{u(W,V)}{v(V,X)})^2 \
   det( (B'[W])^t B'[W])
  \end{array}
 $$
In particular, the coefficients of the above sum over k reduced spanning coforest yielding $\sigma_{b_d - k}( M\{b_r\})$
are all integers.

\end{thm}

Note that if $H_d(X;Q)=0$, then $Z_d(X;Q) = B_d(X;Q)$ and each spanning forest $W$ has complement exactly
a spanning coforest $W^c$ and the complement of a k augmented spanning fores
is a k-reduced spanning coforest. In particular, in this situation
theorems \ref{thm1} and \ref{thm1} become identical.

\section{\textbf{Proof of the main theorem for mesh matrix for $d$ boundaries: Theorem \ref{thm2}}}
\vspace{.1in}

Now consider the case  $V$ a spanning coforest, i.e., $k=0$,  the mappings $k_V,k_V^\star$ are rationally isomorphisms.
Hence, the determinant $ det( k_V \cdot  k_V^\star)$ is the square of the order of the
finite  abelian group:
$$
C_d(X_V;Z)/\{ \ image \ by \ k_V \ of \ B_d(X;Z) \ \}
$$
It is claimed that this finite torsion group has order $v(V,X)$
completing the proof for the case $k=0$.

\vspace{.3in}
Quite generally, let $W$ be a k reduced spanning forest, then
$k_W: B_d(X;Z) \rightarrow C_d(X_V;Z)$ is rationally onto,
so the quotient $C_d(X_W;Z)/\{ \ image \ by \ k_W \ of \ B_d(X;Z) \ \}$
is a torsion group. It is claimed that the order of this finite quotient
group
$$
| \ C_d(X_W;Z)/\{ \ image \ by \ k_W \ of \ B_d(X;Z) \ \} \ |
$$
is $v(W,X)$.

The mapping $k_W$ may be reinterpreted [ since $X_W,X_{W^c}$ each have the same $d-1$ skeleton as $X$:
$$
k_W : B_d(X;Z) \stackrel{onto}{\rightarrow}  B_d(X,X_{W^c};Z) \subset Z_d(X,X_{W^c};Z) = C_d(X_W;Z)
$$
So this finite abelian quotient is identified with $H_d(X,X_{W^c};Z)$. Its order is $v(W,X)$.

Finally, in view of the long exact sequence of the triple $(X,X_{W^c}, X_{V^c})$, there
is an exact sequence
$$
H_d(X_{W^c},X_{V^c};Z) \rightarrow H_d(X,X_{V^c};Z) \rightarrow H_d(X,X_{W^c};Z) \rightarrow H_{d-1}(X_{W^c},X_{V^c};Z)=0
$$
since $X_{W^c}, X_{V^c}$ have the same (d-1) skeleton. By the above all these are finite torsion groups.
Hence, by exactness,
$
u(W,V)/(v(V,X)/v(W,X)) = \frac{u(W,V)}{v(V,X)} \cdot v(W,X)
$
is an integer as desired.

 \vspace{.3in}
 Next consider again as above $V$ a spanning coforest with $W = V-I$
 for $|I|=k$ oriented cells of $V$. That is, $W$ is a k reduced spanning coforest.
 Now the mapping
 $$
 y_{W}: B_d(X_{W^c};Z)= B_d(X_{V^c \cup I};Z)  \rightarrow C_d(X_I ;Z)= C_d(X_{V^c \cup I},X_{V^c};Z)
 $$
 is rationally an isomorphism. Recall that $W^c = V^c \cup I$. It can be factored as
 $$
 B_d(X_{W^c};Z) \stackrel{onto}{\rightarrow} B_d(X_{W^c},X_{V^c};Z)
 \subset Z_d(X_{W^c},X_{V^c};Z) = C_d(X_{W^c},X_{V^c};Z) \ .
 $$
 Hence, the quotient $ C_d(X_I ;Z)/y_{W}(B_d(X_{W^c};Z))$ is a finite
 torsion group and is isomorphic to $H_d(X_{W^c}, X_{V^c};Z)$. Its order
 is  $u(W,V)$.

Let $A'$ be the $n_d \times b_d$ matrix represent the mapping $k_X : B_d(X;Z) \rightarrow C_d(X;Z)$
with respect to the chosen basis $\{b_s\}$ for $B_d(X;Z)$
and the cellular basis for $C_d(X;Z)$. In these terms the boundary mesh matrix,
$M( \{b_s\})$, is just $(A')^t \cdot (A')$.

\vspace{.3in}
Now to evaluate the covolume $det(k_W \cdot k_W^\star)$ for $W$ a k reduced spanning coforest
with $0 < k < b_d$, pick a spanning coforest $V$ of $X$  with $W \subset V$. Let $W = V - I$
with $|I| =k$.

Then by $V$ a spanning coforest the restriction mapping
$$
k_V : B_d(X;Z) \rightarrow C_d(X_V;Z)
$$
is rationally an isomorphism. Also by above, the quotient is of finite order $v(V,X)$.

In particular,   $k_V$ is rationally an isomorphism. so for each $\sigma_j \in V$
[oriented as in $S_d$] there is a unique d boundary  $b(\sigma_j) \in B_d(X;Q)$
which restricts to $1 \cdot \sigma$ in $C_d(X_W;Q)$. That is,
$b(\sigma_j)$ lies in $B_d(X_{V^c \cup \sigma_j};Q) $
Here $\sigma_j$ ranges over $V$ which is of size $b_d$, so $1 \le j \le b_d$.

Under the identification above of $B_d(X;Q)$ with $ C_d(X_V)$,
identifying $b(\sigma_j)$ with $1 \cdot \sigma_j$ for $\sigma_j \in V$,
the map $k_V$ is identified with the identity mapping
$$
Id : B_d(X;Q) \rightarrow B_d(X;Q)
$$
with $A'$ recording the change of basis mapping
$$
b_s = \Sigma_{\sigma_j \in V, j=1, \cdots b_d}\  A'_{\sigma_j, s} b(\sigma_j) \ .
$$

Note that there ia a direct sum decomposition
$$
B_d(X;Q) =  B_d(X_{V^c \cup W},Q)) \oplus ( B_d(X;Q) \cap C_d(X_{V^c \cap I}; Q)
$$
which is compatible with the orthonormal basis given by the geometric generators
$\{ b(\sigma_j) \ | \ \sigma_j \in V\}$. The first summand has basis $\{b(\sigma_j) \ j=1 \cdots b_d-k\}
= \{ b(\sigma_j)\ | \ \sigma_j \in W\}$; the second summand  has basis $\{b(\sigma_j) \ | \ \sigma_j \in I\}$.

Arrange that the first $b_d-k$ cells of $V$ be those of $W = V -I$,
and let $R'$ be the matrix obtained from $A'$ by deleting the last $k$
rows. Then $R'$ represents $k_W$
and $ k_W \cdot k_W^\star$ is represented by $R' (R')^t$.

As before let $S$  be the last $k$ rows of $((A')^{-1})^t$.
Then $S^t$ represents the mapping
$$
B_d(X_{W^c};Z) =  B_d(X_{V^c \cap I}; Q) \subset B_d(X;Q)
$$
with the geometric basis used as a basis for the domain
and the the fixed chosen basis $\{b_s\}$ for the range.

Converting from the geometric basis on the domain
to any choice of integral basis for
$
B_d(X_{W^c},;Z)
$
entails a change with determinant the same as the order
of
$$
C_d(X_I;Z)/y_{W}( B_d(X_{W^c};Z))
$$
which is by definition $u(W,V)$. That is,
$$
det(S S^t) = u(W,V)^2 \ det(B[W]' (B[W]')^t)
$$

Appealing to lemma \ref{lemma55} gives the desired result
$$
\begin{array}{l}
det( k_W \cdot k_W^\star) = det(R' \cdot (R')^t) \\
= det(A')^2 \ det( (S') (S')^t) = det(A')^2 \ u(W,V)^{-2} \ det(B'[W] (B[W]')^t) \\
= (u(W,V)/v(V,X))^2 \ det(B'[W] (B[W]')^t)
\end{array}
$$
as desired.

Since the determinant $det( k_W \cdot k_W^\star) $
as a covolume of an integral lattice in the standard lattice
with standard inner product is an integer, the quantity
$(u(W,V)/v(V,X))^2 \ det(B'[W] (B[W]')^t)$ is an integer also.

The equality $det( k_W \cdot k_W^\star) = (u(W,V)/v(V,X))^2 \ det(B'[W] (B[W]')^t)$
proves that the quotient $f(W) = \frac{u(W,V)}{v(V,X)}$ is independent of  the choice of $V$.

\section{\textbf{Evaluation of all coefficients of mesh matrices for cycles and boundaries
 with respect to   geometrically defined bases.}}
 \label{sectiongeombases}
 \vspace{.1in}

Let $V_0 \subset S_d$ be a fixed choice of d dimensional spanning forest of $X$.
As explained in \S \ref{sectioncyc},
there is a natural geometric basis for the rational $d-$cycles
$Z_d(X;Q)$ described as follows: For each oriented d cell, say $\sigma_j $
in $S_d -V_0$, let $z(\sigma_j)$ be the unique d cycle in $Z_d(F \cup \sigma_j; Q)$
which evaluates to $+1$ on $\sigma_j$. As explained there, as $\sigma_j$ ranges over
$S_d-V_0$  the d cycles
$z(\sigma_j)$ give a basis for the rational vector space $Z_d(X;Q)$
of d cycles on $X$. Also, for any subset $U \subset S_d - V_0$, the
elements $\{ z(\sigma_j) \ | \ \sigma_j \in U\}$ forms a rational
basis for $Z_d(V_0 \cup U; Q)$.

Associated to this geometric basis $\{z(\sigma_j)\}$ one defines the cycle mesh matrix of $X$
for this geometric basis using  $V_0$ as
$M_{cycle, (V_0,X)}\{ z(\sigma_j)\}$  via
$(\sigma,\sigma') \mapsto < z(\sigma),z(\sigma')>_X$ with the pairing
defined using the cell complex $X$.

\begin{thm} \label{thm33} Let $V_0$ be a d dimensional spanning forest, then the cycle
mesh matrix $M_{cycle, (L,X)}\{ z(\sigma_j)\}$ is a $z_d \times z_d$
matrix with coefficient of $(-1)^{k} t^{z_d -k}$ the $k^{th}$ elementary
symmetric function of its eigenvalues, called here $\sigma_k(V_0)$. Moreover,
there is a direct combinatorial formula for this coefficient:
$$
\begin{array}{l}
\sigma_k(V_0)
= \Sigma_{U \subset S_d - V_0 \ with \ |U| =k} \ \  ( \Sigma_{V, \ a \ spanning \ forest \ in \ X_{V_0}} \
( \ \frac{t_{d-1}(X_V)}{t_{d-1}(X_{V_0 \cup U})} \  )^2
\end{array}
$$
\end{thm}

\vspace{.1in}

Analogously, let  $V_1 \subset S_{d+1}$ be a fixed choice of $d+1$ dimensional spanning forest of $X$.
Again, there is a natural geometric basis for the $d$ boundaries $B_d(X;Q)$.
It consists of the boundaries $\{ \partial_{d+1}( 1 \cdot u ) \ | \ u \in S_{d+1}- V_1\}$.
Also, for each subcomplex of the form $X^{(d) } \cup V_1 \cup U$ with $U \subset S_{d+1} - V_1$,
the images $\{\partial_{d+1}( 1 \cdot u) \ | \ u \in S_{d+1}- V_1 \  and  \ u \in U\}$
form a rational basis for $B_d( X^{(d)} \cup V_1 \cup U;Q)$, the boundaries.

Associated to this geometric basis $\{\partial_{d+1}(u) \ | \ u \in S_{d+1}-V_1 \}$,  one defines the boundary  mesh matrix of $X$
for this geometric basis defined by $V_1$ as \newline
$M_{boundary , (V_1,X)}(\{\partial_{d+1}(u) \ | \ u \in S_{d+1}-V_1\})$  via
$(u,v) \mapsto < \partial_{d+1}(u),\partial_{d+1}(v)>_X$ with the pairing
defined using the cell complex $X$.

\begin{thm}  \label{thm34} Let $V_1$ be a $d+1$  spanning forest, then the boundary
mesh matrix  \newline $M_{boundary, (V_1,X)}(\{\partial_{d+1}(u) \ | \ u \in S_{d+1}-V_1\})$  is a $b_d \times b_d$
matrix with coefficient of $(-1)^{k} t^{b_d -k}$ the $k^{th}$ elementary
symmetric function of its eigenvalues, called here $\sigma_{k}(V_1)$. Moreover,
there is a direct combinatorial formula for this coefficient:
$$
\begin{array}{l}
\sigma_{k}(V_1)
= \Sigma_{U \subset V_1 \ with \ |U| =k} \ \ ( \ \Sigma_{V, \ a \ d\ dimensional \ spanning \ coforest \ in \
 X^{(d)} \cup V_1 \cup U}\ |H_{d-1}( X_{V_1}, X_{W^c};Z)|^2
   \ )
\end{array}
$$
where  $ |H_{d-1}( X_{V_1}, X_{W^c};Z)|$ is the order of the finite torsion group  $H_{d-1}( X_{V_1}, X_{W^c};Z)$.

with $X_{V_1} = X^{(d)} \cup V_1$

and $X_{W^c} = X^{(d-1)} \cup (S_d- W)$.
    This weight is alternately expressed  as

   the integer $ [ \  \frac{ t_{d-1}(X_W) t_{d-2}(X_{(V_1)^c}) t'_{d-1}(V_1)}{ t_{d-2}(X)} \ ]$,
   see the next section.
\end{thm}

The proofs of these two theorems  proceeds as follows:

For theorem \ref{thm33}, applying lemmas \ref{lemmaalg1}, \ref{lemmaalg2},
one sees that $\sigma_k(V_0)$ is expressed as a double sum over $U \subset S_d - V_0$ with $|U| =k$
and over $Y \subset S_d$ of size $k$ of a square of the determinant for the restriction mapping
$Z_d(V_0 \cup U;Q)$ to $C_d(Y;Q)$. This last mapping is onto [and hence an isomorphism]
only  if $Y \subset (V_1 \cup U)$ and $V = (V_1 -Y)$ is a spanning forest for $X_{V_1 \cup U}$.
In  this instance, the matrix at issue is that for changing from the geometrical basis
for the spanning forest $V_0$ in $X_{V_0 \cup U}$  to the geometrical basis for the spanning forest  $V$
of the same cell complex  $X_{V_0 \cup U}$. These two bases can be compared to the integral basis of $Z_d(X_{V_0 \cup U};Z)$.
by invoking  theorem \ref{thm1}. On that basis, the desired determinant is
$(\frac{t_{d-1}(X_{V_0})}{t_{d-1}(X_{V_0 \cup U})})^{-1} \ (\frac{t_{d-1}(X_{V})}{t_{d-1}(X_{V_0 \cup U})})
= \frac{t_{d-1}(X_V)}{t_{d-1}(X_{V_0})}$.

For theorem \ref{thm34}, applying lemmas \ref{lemmaalg1}, \ref{lemmaalg2},
one sees that $\sigma_k(V_1)$ is expressed as a double sum over $U \subset V_1$ with $|U| =k$
and over $V \subset S_d$ of size $k$ of a square of the determinant for the  mapping
$$
C_{d+1}(X^{(d)} \cup V_1,Z) \stackrel{\partial_{d+1}}{\rightarrow }
   B_d(X^{(d)} \cup V,Z) \subset C_d(X;Z) \stackrel{restriction}{\rightarrow}  C_d(X^{(d)} \cup V,Z) \ .
   $$

 This last mapping is onto [and hence an isomorphism]
only  if $V \subset X^{(d)} \cup U$ and $V$ is a d dimensional spanning coforest for $X^{(d)} \cup U$.
As the next section proves, in this case the mapping
   is rationally an isomorphism with quotient group  a  finite order equal to the order of the torsion group
   $H_{d-1}( X_{V_1}, X_{W^c};Z)$, Here  $X_{V_1} = X^{(d)} \cup V_1$ and $X_{W^c} = X^{(d-1)} \cup (S_d- W)$.
    This is alternately expressed there as
   the integer $ [ \  \frac{ t_{d-1}(X_W) t_{d-2}(X_{(V_1)^c}) t'_{d-1}(V_1)}{ t_{d-2}(X)} \ ]$.

\section{\textbf{A Kirchhoff type  theorem for cell  complexes, extending  Lyons's \cite{Lyons}   work
to the complete characteristic polynomial }}  \label{sectionKirk}
\vspace{.1in}

In the same manner:

Let $A$ be the integer valued matrix representing  the boundary mapping  with respect to the cellular basis
$$
\partial_d : C_d(X;Z ) \rightarrow C_{d-1}(X;R).
$$
Then $A^t$ represents the adjoint for the cellular inner product $<.,>$. Recall that
$\partial_d $ has  rank $b_{d-1} = dim \ B_{d-1}(X;R)$

\vspace{.1in}
Let $0 < m \le b_{d-1}$ with $S_d$ the $d$ cells with their specified
orientations.

By the same algebra as before, one gets
$$
\begin{array}{l}
m^{th} \ elementary \ symmetric \ function \ of \ the \ eigenvalues \ of \ \partial \cdot \partial^\star \\
= \Sigma_{I \subset S_{d-1}\ J \subset S_d \ with \ |I|=|J| = m} \ \ [ \  det( A_{(I,J)}\ ]^2 \ .
\end{array}
$$
For the  term to be non-zero, the mapping $A_{(I,J) }$ must be an isomorphism after tensoring
with the reals.

The matrix $A_{(I,J)}$ represents the composite :
$$
C_d(X_J;Z) \subset C_d(X;Z) \stackrel{\partial}{\rightarrow} B_{d-1}(X;Z) \subset C_{d-1}(X;Z) \stackrel{restrict}{\rightarrow} C_{d-1}(X_I;Z)
$$

For $det(A_{(I,J)}) \neq 0$, necessarily the composite
$$
  C_d(X_J;R) \subset C_d(X;R) \stackrel{\partial}{\rightarrow} B_{d-1}(X;R)
  $$
must be one to one. That is, $Z_d(X_J;Z)=0$.

The conventional notation
for this is : $J \subset S_d$ is called a \textbf{forest of size $|J|=m$}.
That is, $X_J$ carries no $d$-cycles.

Additionally, for  $det(A_{(I,J)}) \neq 0$, necessarily the composite
$$
B_{d-1}(X_J;R) \subset   B_{d-1}(X;R) \subset C_{d-1}(X;R) \stackrel{restrict}{\rightarrow} C_{d-1}(X_I;R)
  $$
must be an isomorphism. That is, $I\subset S_{d-1}$ interpreted via $X_I = X^{d-2} \cup I$ inside $X_J = X^{d-1} \cup J$
 must be  a spanning coforest of $X_J$. In combination this implies that $|I|= m$ also.

\vspace{.1in}

 Hence, the above sum is restricted to $J=V$ a (d-dimensional) forest of size $m$  and $I=W$ a
 ((d-1)-dimensional)  spanning coforest of $X_J= X_V$.
 This ensures that $det(A_{(I,J)})  \neq 0$.
 $$
\begin{array}{l}
m^{th} \ elementary \ symmetric \ function \ of \ the \ eigenvalue \ of \ \partial \cdot \partial^\star \\
= \Sigma_{V \subset S_{d}\ W \subset S_{d-1} \ with \ |V|=m,   \ |W| =m ,  \ V \ a \ forest\ of \ size \ m \ and \  W \ a
\ spanning \ coforest \ of \ X_V
 } \ \ [ \  det( A_{(W,V)}) \ ]^2
\end{array}
$$

Consider the inclusion $X_W = X^{d-2} \cup W \subset X_V = X^{d-1} \cup V$, defining the
mapping $A_{(W,V)}$
$$
C_d(X_V;Z) \subset C_d(X;Z) \stackrel{\partial}{\rightarrow} B_{d-1}(X;Z) \subset C_{d-1}(X;Z) \stackrel{restrict}{\rightarrow} C_{d-1}(X_W;Z)
$$
which in the cases of interest is a rational isomorphism.

Now as A.  Duval, C. Klivans,  and J.  Martin note \cite{DuvalKlivansMartin},
one can reinterpret   this map more geometrically. On integral chains,
$$
C_d(X_V,X_{W^c};Z) = C_d(X_V;Z) \stackrel{\partial }{\rightarrow} C_{d-1}(X_V) \stackrel{restrict}{\rightarrow}
 C_{d-1}(X_W;Z) = C_{d-1}(X_V,X_{W^c};Z)
 $$
  becomes identified with the boundary mapping for the pair $(X_V,X_{W^c})$. Since $C_{d-1}((X_V,X_{W^c};Z) =0$
 and $C_{d+1}((X_V,X_{W^c};Z)=0$ and this boundary mapping is a rational isomorphism, it follows that
 $$
 det(A_{(V,W)})^2  = |H_{d-1}(X_V,X_{W^c};Z)|^2
 $$
 in the cases of interest.

\vspace{.1in}
In the maximal case $k = b_d$, by $V$ a forest of size $b_d$, necessarily the mapping
$C_d(X_V;Z) \subset C_d(X;Z) \stackrel{\partial}{\rightarrow} B_{d-1}(X;Z)
$ is a rational isomorphism. Hence, the set of $W$ for which $W$ is a spanning
coforest of $X_V$ is the same as the spanning coforests  of $X$.

Lyons proved in the  case $k = b_d$, that
$$
|H_{d-1}(X_V,X_{W^c};Z)|^2 = [ \  \frac{ t_{d-1}(X_W) t_{d-2}(X_{V^c}) t'_{d-1}(V)}{ t_{d-2}(X)} \ ]^2 \in Z
$$
where for a CW complex $K$ the order of the torsion of $H_k(K;Z)$ is denoted by $t_k(K)$
and $t'_{d-1}(V) $ is the order of the quotient of $Z_{d-1}(X;Z) $
by  $Z_{d-1}(X) \cap B_{d-1}(X;Q) + Z_d(X_ {V^c};Z)$. Here by convention
$X_{W} = X^{(d-1)} \cup W$ and $X_{V^c} = X^{(d)} \cup V^c$.

For $m=b_d$ the corresponding sum decouples, yielding the first part of the next theorem,
which appears in the work of Lyons:

\begin{thm}

\begin{itemize}
\item[a.] Lyons \cite{Lyons}.  The product of the  non-zero  eigenvalues  of  $\partial \cdot \partial^\star$,
the combinatorial Laplacian, or Kirchhoff mapping,
acting on (d-1) chains $C_{d-1}(X;Z) \rightarrow C_{d-1}(X;Z)$, i.e., equals
$$
\begin{array}{l}
= (\ \Sigma_{W, \ forest \ of \ size \ m = b_d}  \  (t_{d-1}(X_W))^2 \ ) \cdot
( \ \Sigma_{V, \ spanning \ coforest \ in \ X }  \ ( \frac{  t_{d-2}(X_{V^c}) t'_{d-1}(V)}{ t_{d-2}(X)}  )^2 \ ),
\end{array}
$$
where each term of each sum is an integer.

\item[b.] More generally,  if $0 <  k \le b_d$,  the $k^{th}$ elementary symmetric function
in  the eigenvalues of $\partial \cdot \partial^\star$, the combinatorial Laplacian,
or Kirchhoff map is
$$
\begin{array}{l}
=  \Sigma_{V \ a \  \ forest \ of \ size \ k }   \ [ \
 \Sigma_{W \ a \ spanning \ coforest \ of \ X_V
 }  \ |H_d( X_V, X_{W^c};Z)|^2 \ ] \ .
\end{array}
$$
\end{itemize}
\end{thm}

\section{\textbf{Homology covolume squared via determinants of mesh matrices}}  \label{sectionReid}
\vspace{.1in}

Naturally, the above results imply a statement about $H_d(X;R)$.
Recall that using the cellular basis to identify  chains and cochains,
the coboundary dual to $\partial_d:C_d(X;R) \rightarrow C_{d-1}(X;R)$ gets identified with
 $\partial_d^\star : C_{d-1}(X:Z) \rightarrow C_d(X;Z)$.
The Hodge theory decomposition becomes in this setting  an orthogonal direct sum decomposition
of real $d$ chains with the cellular inner product $<.,>$:
$$
C_d(X;R) = \mathcal{H}_d \oplus \{ image \ \partial_{d+1} \} \oplus \{ image \ \partial_d^\star \},
$$
with $\mathcal{H}_d$ the space of (d-dimensional) harmonics $\mathcal{H}_d = \{ c \ | \ \partial_d(c)=0 \ and \ \partial_{d+1}^\star(c)=0\}$.

Now one has the direct sum orthogonal decomposition :
$$
Z_d(X;R) =  \mathcal{H}_d \oplus B_d(X;R)
$$
Let $\overline{B_d(X;Z)} =  B_d(X;Q) \cap Z_d(X;Z)$. Then $\overline{B_d(X;Z)}$ is a direct summand
of $C_d(X;Z)$ and the quotient $\overline{B_d(X;Z)} / B_d(X;Z)$ is isomorphic to the torsion in $H_d(X;Z)$.

Thus the covolume squared  for the lattice $B_d(X;Z)$ in $C_d(X;Z)$ with cellular inner product $<.,>$
is $t_d(X)^2$ times that of $\overline{B_d(X;Z)}$.

Now chose a direct summand, say $G \subset Z_d(X;Z)$, with
$$
Z_d(X;Z) = G \oplus \overline{B_d(X;Z)} \ .
$$
This is feasible as $\overline{B_d(X;Z)}$ is a direct summand.
Then the covolume of $Z_d(X;Z)$ computed from this lattice is exactly the same
as that of $\pi(G) \oplus \overline{B_d(X;Z)}$, where $\pi_d$ is is the orthogonal projection
of $Z_d(X;Z)$ onto $\mathcal{H}_d$ or equivalently the direct summand $G$.
This shows that
$$
\begin{array}{l}
covolume  \ squared \ for \ lattice \ Z_d(X;Z) \ under \ <.,.> \\
= ( \ covolume  \ squared \ for \ lattice \ B_d(X;Z) \ under <.,.> \ ) \\
\hspace{.5in} \cdot ( \ covolume \ for \ the \ lattice \ \pi_d(Z_d(X;Z)) \subset \mathcal{H}_d  \ under \ <.,>\ ).
\end{array}
$$

That is, if  one   defines the lattice in $\mathcal{H}_d$
by the orthogonal projection $\pi_d$ of the lattice $Z_d(X;Z)$ onto $\mathcal{H}_d$, one has
$$
\begin{array}{l}
covolume  \ squared \ for \ lattice \ \pi_d(Z_d(X;Z)) \subset \mathcal{H}_d  \ under <.,.> \\
= \frac{ covolume  \ squared \ for \ lattice \ Z_d(C;Z) \ under <.,.> }{ covolume  \ squared \ for \ lattice \ B_d(C;Z) \ under <.,.> } \cdot t_d(X)^2 \ .
\end{array}
$$

This gives in view of the above a combinatorial description of this  ``homological covolume''.

\section{\textbf{Combinatorial identities via a computation of Reidemeister-Franz torsion}} \label{sectRF}
\vspace{.1in}

As explicated  by Milnor \cite{Milnor2},
the R-torsion of Reidemeister-Franz \cite{Reidemeister,Franz} \newline
$\tau(C_\star, \{d_i\},\{e_{i,j}\}, \{h_{j,k}\})$,
 is an invariant which is a non-zero real number, associated to the following data:
\begin{itemize}
\item A finite chain complex of real finite dimensional vectors spaces, say $C_0,\cdots, C_N$
with differentials $d_i : C_i \rightarrow C_{i-1}$ with $d_{i-1} \cdot d_i = 0$ as recorded by
the sequence:
$$
0   \rightarrow C_n \stackrel{d_n}{\rightarrow} C_{n_1} \stackrel{d_{n-1}}{\rightarrow } \cdots
  \stackrel{d_1}{\rightarrow} C_ 0 \rightarrow 0
  $$
    Let $C_i$ have real dimension $r_i$. Let $Z_i =\{kernel \ of \ d_i :C_i \rightarrow C_{i-1}\}$
  be the i-cycles; let  $B_i = \{ image \ of \ d_{i+1} ;C_{i+1} \rightarrow C_i\}$ be the i-boundaries.
  Let $H_i$ be the $i^{th}$ homology,
  $H_i = Z_i/B_i$
  and have rank $s_i$. Let $B_i$ have rank $t_i$.

  \item For each vector space $C_i$ a chosen basis, say $\{e_{i,j} \ | \ j=1,\cdots, r_i\}$
  \item For each homology $H_i$ a chosen basis, say $\{ h_{i,k} \ | \ k=1,\cdots, s_i\}$
  \end{itemize}

  This torsion $\tau$ is defined as follows: Chose elements $\tilde{h}_{i,k} \in Z_i \subset C_i$
  which project to the chosen homology classes $h_{i,k}$. Next chose elements, say $ \{ b_{i-1,l} \ | \ l=1,\cdots, t_{i-1}\}$
  which are a basis for $B_{i-1}$. Next chose elements $\tilde{c}_{i,l} \in C_i$ with $d_i(\tilde{c}_{i,l}) = b_{i-1,l}$
  with $l=1,\cdots, t_{i-1}$.

  With these choices, $Z_i$ has basis $[\{ \tilde{h}_{i,k}\} \sqcup  \{ b_{i,l} \}]$ and
  $C_i$ has the basis \newline
   $X_i = [ \{ \tilde{h}_{i,k}\} \sqcup  \{ b_{i,l} \} \sqcup \{ \tilde{c}_{i,l'} \}]$.
  Now set $D_i = |det( [ X_i: \{e_{i,j}\}])|$ where $[ X_i: \{e_{i,j}\}]$ is the matrix expressing the basis $X_i$
  in terms of the basis $\{e_{i,j}\}$.
   One defines the Reidemeister-Franz   torsion via :
  $$
  (\tau(C_\star, \{d_i\},\{e_{i,j}\}, \{h_{j,k}\}))^2 = \prod_{i=0}^N \ (D_i)^{(-1)^i}
  $$
   A simple check shows that the R-torsion $\tau$ squared is independent of the choices of
  $\tilde{h}_{i,k}, b_{i-1,l}, c_{i,l'}$ \cite{RaySinger} ; also,
  the sign issues are  disposed of   by taking the square.

  \vspace{.1in}
  Turning to the case of a cell complex of dimension $d$, one has the cellular basis for $C_i(X;R)$  given
  by the i-cells of $X$, once oriented. Moreover, the real homology $H_i(X;R)$ has a preferred basis
  coming from a choice of $Z$-basis, say $\{h_{i,k}\}$ for the free part of the integral homology.
  Consequently,  there is a canonically defined R-torsion for any finite cell complex $X$,
  $\tau( C_\star(X;R), \partial_R, \{ combinatorial \ basis \ for \ C_\star(X;R) \ and \ H_\star(X;R)\})$.
  Not surprisingly, there is a simple formula (known to experts) for it:

\begin{thm} \label{thm4}
The combinatorial R-torsion is given by
$$
[ \ \tau( C_\star(X;R), \partial_R, \{ combinatorial \ basis \ for \ C_\star \ and \ H_\star) \} \ ]^2
    = [ \ \prod_{k=1}^{d-1}  \theta_{X,i}^{(-1)^i} \ ]^2 \ ,
    $$
    where $\theta_{X,i}$ is the order of the torsion in $H_i(X;Z)$.
    \end{thm}

    Proof: Since $Z_i(X;Z)$ is a free direct summand of $C_i(X;Z)$, the quotient \newline $C_i(X;Z)/ B_i(X;Z)$
    has torsion subgroup equal to the torsion subgroup  of $H_i(X;Z) = Z_i(X;Z)/B_i(X;Z)$;  so in particular
    its order  is just $\theta_{X,i}$. Alternatively expressed, let $\overline{B_i(X;Z)}$ be the lattice
    of elements $x \in C_i(X;Z)$ with  $N c \in B_i(X;Z)$ for an integer $N \neq 0$, then
    $ |\overline{B_i(X;Z)}/B_i(X;Z) | =\theta_{X,i}$ and $\overline{B_i(X;Z)}$ is a direct summand of $C_i(X;Z)$.

     Now
    chose elements say $ \hat{b}_{i,l} \in \overline{B_i(X;Z)}$ that are a  $Z$ basis for the direct summand
   $\overline{B_i(X;Z)}$.
    Since $Z_i(X;Z)/\overline{B_i(X;Z)} \cong H_i(X;Z)/\{ torsion \} $ is free abelian, one may choose classes
    $ \{ \tilde{h}_{i,k} \}$ in $Z_i(X;Z)$ which project to a $Z$-basis for $H_i(X;Z)/\{torsion\}$, i.e.,
    the required combinatorial basis. With these choices made
    $\{ \tilde{h}_{i,k} \} \sqcup \{ \hat{b}_{i,l} \}$ form an integral basis for
    the direct summand $Z$-lattice $Z_i(X;Z)$ inside $C_i(X;Z)$. Again, since $C_i(X;Z)/Z_i(X;Z)$ is free abelian,
    one may chose elements $x_{i,l'} \in C_i(X;Z)$
    so that they project to an integral basis of this quotient.
    With these choices the  lattice $C_i(X;Z)$ has an  integral basis given by
     $\{ \tilde{h}_{i,k}\} \sqcup \{\hat{b}_{i,l}\} \sqcup \{ x_{i,l}\}$. This gives a reference basis
    for $C_i(X;Z)$ and represents in $C_i(X;R)$ a basis which differs from the cellular basis
    by an element in $GL(dim \ C_i(X;R),Z)$ so the transition matrix has determinant $\pm 1$.
    Consequently, one may use this reference basis in computing the $D_i$ and the R-torsion $\tau$.

    Naturally, one  uses the chosen lift $\tilde{h}_{i,k}$ again and the   classes
     $\{ d_{i+1}(x_{i+1,l'} )\}$ and $\{x_{i,l}\}$  in computing the R-torsion with  the above
    reference basis. That is, \newline $X_i = [\{ \tilde{h}_{i,k} \}, \{ d_{i+1}(x_{i+1,l'}) \}, \{ x_{i,l} \}]$
    is a basis for $C_i(X;R)$ which conforms to the above recipe. $X_i$ is a sublattiice
    of the reference lattice $ [\{ \tilde{h}_{i,k}\}, \{ \hat{b}_{i,l'}\},\{ x_{i,l}\}] $
    of index $|\overline{B_i(X;Z)}/B_i(X;Z)| =\theta_{X,i}$.
    Hence, $D_i = \theta_{X,i}$ and
    the theorem follows.

\vspace{.1in}
Returning to the general context of the R-torsion, a chain complex with the above choices $\{x_{i,j}\}, \{ h_{i,k}\}$,
 there is a real inner product
in which these basis elements are orthonormal. Hence, for these inner products the boundary mapping
$d_{i+1}$ has an adjoint
$$
d_{i+1}^\star : C_i \rightarrow C_{i+1} \ .
$$
As before let $\{h_{i,k}\} $ be a specified basis for the real vector space $C_i$.

As Ray and Singer remark in their foundational paper on analytic torsion
\cite{RaySinger}, there is a formula for the
R-torsion $\tau$ in terms of the reduced Laplacian
$$
\tilde{\Delta}_{X,i}  := (d_{i+1} \cdot  d_i^\star|B_i) : B_i \rightarrow B_i \ .
$$

The first part of this relation is a measure, via covolumes
of the chosen homology basis:
For this purpose, note that one has an orthogonal Hodge decomposition,
$$
C_i = Harm_i \oplus \{image \ of \ d_{i+1}\} \oplus \{ image \ of \ d_{i-1}^\star \} \ ,
$$
in which $Harm_i = \{ h \ | \ d_i(h) =0 \ and \ d_{i-1}^\star(h) =0 \}$ are the i-harmonics.
Moreover, $Z_i =Harm_i \oplus \{ image \ of \ d\}$ is a orthogonal decomposition.
Let $\pi_i$ be the orthogonal projection :
$$
\pi_i: C_i  \rightarrow Harm_i
$$
It induces  the  isomorphism $H_i(X;R) =  Z_i/B_i \stackrel{\cong}{\rightarrow } Harm_i$\ .

In this manner the specified basis for the integral homology modulo the torsion $\{h_{i,k}\}$
under the orthogonal projection  to $Harm_i$
yields  the desired elements
 $\pi_i(h_{i,k}) \in Harm_i$. They define a $Z$-lattice of full rank
 which are independent of the lifts chosen for the $\{h_{i,k}\}$'s. The subspace $Harm_i$ inherits
an inner product from that on $C_i$. One defines
$$
(\mathcal{H}_{X,i})^2 = [covolume( lattice \ \{\pi_i(h_{i,k})\} \ in \ C_i)]^2
= det( \{ < \pi_i(h_{i,a}), \pi_i(h_{i,b})> \} )
$$
as a measure of the size of the combinatorial homology basis.

With these preparations one has:

\begin{thm}[Ray and Singer \cite{RaySinger}]

The R-torsion is given in terms of the covolumes $\mathcal{H}_i$ for
the i-homologies and the determinant of the reduced Laplacians by
$$
\begin{array}{l}
[ \ \tau( C_\star(X;R), \partial_R, \{ combinatorial \ basis \ for \ C_\star(X;R) \ and \ H_\star(X;R)) \ ]^2 \\
= \prod_{i=0}^n  \ [  det( \tilde{\Delta}_{X,i}) \cdot (\mathcal{H}_{X,i})^2] ^{(-1)^i} \ .
\end{array}
$$
\end{thm}

[Note that in Ray and Singer's paper they use the full Laplacian $\Delta = d_{i+1} \ d^\star_i + d_{i-1}^\star \ d_i$
with corresponding modification of the above  formula. They also restrict to the acyclic cases
$H_i =0$.]

\vspace{.1in}
Returning to the cell context and the combinatorial bases for chain and homologies,
The formula of Ray and Singer and the above computation yields:
\begin{thm} For each finite $CW$ complex, and each choice of skeleton, say $X^{(d)}$,
$$
\begin{array}{l}
 [ \ \prod_{k=1}^{d-1}  \theta_{X^{(d)},i}^{(-1)^i} \ ]^2\\
= \prod_{i=0}^d  \ [  det( \partial_{i+1} \cdot \partial_i^\star | image \ of \ \partial_{i+1}) \cdot (\mathcal{H}_{X,i})^2] ^{(-1)^i} \ .
\end{array}
$$
\end{thm}

This theorem gives a sequence of combinatorial identities
since the terms $det( \partial_{i+1} \cdot \partial_i^\star|image \ of \ \partial_{i+1})$
are expressed by Lyons's formulae combinatorially.

\vspace{.1in}
The papers of  M. Catanzaro, V. Chernyak,
and  J. Klein, \cite{CatanzaroChernyakKlein1, CatanzaroChernyakKlein1,CatanzaroChernyakKlein2}
also study aspects of the Reidemeister-Fraz torion in the combinatorial
context. Their results are complementary to those given here.

\section{\textbf{Variant Settings and Examples}}
\vspace{.1in}

\subsection{Weights}

One may easily adapt the above discussion to the context of cell complexes with
weighted cells. Let $\{e_{i,j}\}$ be the i-dimensional cells of a cell complex $X$
and $W_i : \{e_{i,j}\} \rightarrow R^+$ be a positive  real valued function.
This weighting  defines using the cellular basis an isomorphism, $W_i : C_i(X;R) \rightarrow C_i(X;R)$
multiplying a cell by its weight.
Some familiar  examples for graphs are:
\begin{itemize}
\item  One often  introduces multiplicities for the edges.
\item For the normalized
reduced graph Laplacian a vertex $v$ is  weighted by one divided by the degree of the vertex, see F. R. K. Chung
 \cite{Chung}.
 \item In the classical
Kirchoff circuit setting, the resistances of the wire edges are used as the weights.
\end{itemize}
The paper of M. Catanzaro, V. Chernyak,
and  J. Klein \cite{CatanzaroChernyakKlein2}   generalizes classical results of Kirchhoff on spanning trees, see also
\cite{Bryant}, to  circuits, voltages, resistances, and is written in a weighted higher dimensional setting.

Now  all the present results readily extend to the weighted context. There are two formally
equivalent methods  to carry this out.

In the first method: keeping the
 the standard cellular inner product, one replaces the
mapping $\partial : C_{i+1}(X;R) \rightarrow C_i(X;R)$ represented by the matrix $A$
and its adjoint $\partial^\star : C_i(X;R) \rightarrow C_{i+1}(X;R)$ represented by
the matrix $A$ by the mappings $ \mathcal{L}_W = (W_i)^{-1/2} \partial (W_{i+1})^{1/2}$ represented
by $(W_i)^{-1/2}\cdot A \cdot (W_{i+1})^{1/2}$ and its adjoint  $\mathcal{L}^\star
= (W_{i+1}^{1/2} \cdot \partial^\star \cdot (W_i)^{-1/2}$ represented by $ (W_{i+1}^{1/2} \cdot A \cdot (W_i)^{-1/2}$.
That is, one replaces the reduced Laplacian $\partial \cdot \partial^\star$ represented by $A \cdot A^t$
by $ \mathcal{L} \cdot \mathcal{L}^\star$ represented by
$ ((W_i)^{-1/2}\cdot A \cdot (W_{i+1})^{1/2}) \cdot ( (W_{i+1})^{1/2}\cdot A^t \cdot (W_{i})^{-1/2})
= (W_i)^{-1/2} \cdot A \cdot W_{i+1} \cdot A^t \cdot (W_i)^{-1/2}$.

In this approach the operators are changed, but the inner product is left alone.
All covolumes are computed for the new operators with the old inner product.
With these slight changes, all theorems can be recast in this extended context.

In the second  alternative approach: One    recognizes that the characteristic polynomial
of $ (W_i)^{-1/2} \cdot A \cdot W_{i+1} \cdot A^t \cdot (W_i)^{-1/2}$ equals
that of $ A \cdot W_{i+1} \cdot A^t \cdot (W_i)^{-1}$. Now if takes the adjoint
of $\partial$ represented by $A$ with respect to the modified inner product
$<..>_W$ given by $<e_{i,j},e_{i,k}>_W = \delta_{j,k} \ W(e_{i,j})$
then the adjoint, called here $\partial^{\star,W}$ is represented by
$W_{i+1} \cdot A^t \cdot (W_i)^{-1}$. That is, one may rephrase the
above theorem in terms of $\partial \cdot \partial^{\star,W}$
by computing the covolumes utilizing the new inner product
$<.,>_W$ but with $\partial $  unchanged.

All the results of the present paper can be readily extended to these  weighted settings.
For example, the analogue  of Theorem \ref{thm1} is

\begin{thm} \label{thm15}
Let $X$ be a $d$-dimensional cell complex.

 The $k^{th}$ elementary symmetric polynomial in the eigenvalues of the reduced weighted Laplacian
  $\mathcal{L} \cdot \mathcal{L}^\star$ are combinatorially expressed
  in term of covolumes of lattices  with respect to the weighted metric $<.,>>_W$ as follows:
  $$
  \begin{array}{l}
\sigma_k( \tilde{\Delta}) = [covolume( lattice \ image  \ \partial_R^\star(C_d(X;Z)) \subset C_{d+1}(X;Z)  \ in \  C_{d+1}(X;R))]^2 \\

  =\Sigma_{ forests \ F  \ of \ X \ of \ size \ k } \ \
 ( \Sigma_{coforests \ K \ of \ size \ k \ in \ forest \ F } \\
  \hspace{.5in}   [covolume(  lattice \ image \  B_{d-1}(F;Z) \rightarrow  C_{d-1}(K;R))]^2 ) \\
  \end{array}
  $$
  where the mapping $ B_{d-1}(F;Z) \rightarrow C_{d-1}(K;R)$ proceeds by inclusion into $C_{d-1}(F;R)$ and then restriction to $C_{d-1}(K;R)$

Here if there are no forests, or coforests of the required type, then that term is to be
interpreted as zero.
\end{thm}

And the  analogue of the corollary of the Trent type theorem in the  weighted setting is:
Let all covolumes in the next statement  be computed with respect to the  weighted inner product $<.,.>_W$.
 Then an inductive calculation of the covolumes of $Z_i(X;Z) \subset C_i(X;R)$ computed
 with respect to the weighted inner product is:

\begin{corollary}
$$
\begin{array}{l}
 [covolume \ of \ lattice  \ Z_d(X;Z) \subset C_d(X;R))]^2 \\
 = [covolume \ of \ lattice  \ Z_{d-1}(X;Z) \subset C_{d-1} (X;R))]^2 \ \cdot \\
 \hspace{.3in} ((\mathcal{H}_{k,X})^2 \cdot [covolume( \overline{B_k(X;Z)} \subset C_k(X;R))]^2
                     \cdot (  \Sigma_{T \in \mathcal{T}_{k+1}} \ (\frac{\theta_{T,k}}{\theta_{X,k}})^2) \ .
                     \end{array}
                     $$
\end{corollary}

We leave as an exercise to the reader the rewriting of the other
theorems  above into  this enlarged setting.

Remark: Likewise, all the results of this paper can also be rewritten to encompass
weights which are allowed to be non-zero real numbers (or even non-zero complex numbers),
e.g., using the normal principal values for square roots. An interesting example of
 positive and negative edge weights $\pm 1$ occurs in the work of
K. Marasugi \cite{Murasugi} on polynomial invariants of knots.

\vspace{.1in}
\subsection{Twisted Coefficients}

Historically, Reidemeiser-Franz torsion was introduced as an invariant
of finite cell complexes with $U(1)$ representations of the fundamental
group, equivalently flat $U(1)$ bundes. The definitions were later
extended to arbitrary unitary representations.

From this point of view the above combinatorial results
are in the setting of the trivial $U(1)$ representation.
As also done for some results on spanning trees by
 M. Catanzaro, V. Chernyak,
and  J. Klein \cite{CatanzaroChernyakKlein1},
in this section the above combinatorial results
are extended to $U(1)$ representations.

Given a representation $\rho : \pi_1(X) \rightarrow U(1)$
of the fundamental group of a connected cell complex $X$   into $U(1)$,
one may form the chains with values in $C$ twisted by the
flat representation $\rho$, $C_\star(X;E_\rho)$ together with
it boundary mapping
$$
\partial_\rho : C_{i+1}(X;E_\rho) \rightarrow C_i(X,E_\rho)\ .
$$
Since $\rho$ is unitary, the standard Hermitian pairing
$<z,w> \mapsto z \overline{w}$ on $C$ extends to
a Hermitian inner product on i-chain with values in
the local system  $E_\rho$ :
$$
<.,>: C_i(x;E_\rho) \times C_i(X;E_\rho) \rightarrow C
$$
Hence, the adjoint $\partial_\rho^\star; C_i(X;E_\rho) \rightarrow C_{i+1}(X;E_\rho)$ is well defined.

\vspace{.1in}
Now suppose that $X$ is a d-dimensional cell complex
with i-dimensional cells $\{ e_{i,j}\}$.
One  may consider the one dimensional vector space $C_i(X;E_\rho)$
of chain with support in that cell, say
$$
V[e_{i,j}] = \{ c \in C_i(X;E_\rho) \ | \ c \ vanishes \ on \ cells \ not \ equal \ to \ e_{i,j} \} \ .
$$
This provides a canonical orthogonal direct sum decomposition of
the $i-chains$ :
$$
C_i(X;E_\rho) = \bigoplus_j \  V[e_{i,j}]
$$

For example, for any subcomplex $S \subset X$ containing say $k$ cells of
dimension $d$ and all of the $d-1$ skeleton, one may consider the
vector space of dimension $k$ of $C_d(X;E_\rho)$ define by
$$
V_i[S] = \{ c \in C_i(X;E_\rho) \ | \ c \ vanishes \ on \ cells \ not \ equal \ to \ e_{i,j} \}
    = \bigoplus_{e_{d,j} \ a \ d-cell \ of \ S} \ V[e_{d,j}] \ .
    $$
Naturally there are boundary mapping $(\partial|S): V_{d}[S] \rightarrow V_{d-1}[S]$.

Now picking bases elements of length one for the one dimensional vector
spaces $V[e_{d,j}]$ and $V[e_{d-1,k}]$ the boundary mapping
$$
\partial :C_d(X;E_\rho) \rightarrow C_{d-1}(X;E_\rho)
$$
is represented by a matrix, say $A$, with rows and columns indexed
by the $d-1$-cells of $X$ and the $d-1$ cells of $X$ respectively,
or alternatively the one dimensional blocks, $V[e_{d-1,k}], V[e_{d,j}]]$,
respectively.

Consequently, the algebraic lemmas \ref{lemmaalg1},\ref{lemmaalg2} apply, giving by the usual method

$$
\begin{array}{l}
\sigma_{k}( \partial_\rho \cdot \partial_\rho^\star: C_{d-1}(X;E_\rho) \rightarrow C_{d-1}(X;E_\rho)) \\
 =  \sigma_{k}( \partial_\rho^\star  \cdot \partial_\rho: C_d(X;E_\rho) \rightarrow C_d(X;E_\rho)) \\
 = \Sigma_{Blocks \ of \ size \ k \ V[e_{d,i_1}] \oplus \cdots \oplus V[e_{d,i_k}]} \ \ \
                 det( \partial_\rho^\star  \cdot \partial_\rho | ( V[e_{d,i_1}] \oplus \cdots \oplus V[e_{d,i_k}]) \ .
                 \end{array}
                 $$

\vspace{.1in}
Analogous to the above,
define a subcomplex $F \subset X$ to be a d-dimensional forest of size $k$  for $\rho$
if $F$ includes the $d-1$ skeleton of $X$, $F$ contains exactly $k$ cells of dimension $d$
and the boundary mapping
$$
(\partial|S) : C_d(S;E_\rho) \rightarrow C_{d-1}(S;E_\rho)
$$
is one to one.

By application of the algebraic identity above, one gets the
analogue  of theorem \ref{thm1} in  this twisted context:
$$
\begin{array}{l}
\sigma_{k}( \partial_\rho \cdot \partial_\rho^\star: C_{d-1}(X;E_\rho) \rightarrow C_{d-1}(X;E_\rho)) \\
 = \Sigma_{d- dimensional \ forests \ F\ of \ size \ k} \
                 det( A_F^\star \cdot A_F)
                 \end{array}
                 $$
                 where $A_F $ is the composite
                 $C_d(F;E_\rho) \stackrel{\partial|S}{\rightarrow} C_{d-1}(S;E_\rho ) \subset C_{d-1}(X;E_\rho)$.

                 Of course, there are analogous theorems about forests and coforests and their dual notions.

As mentioned above, the spanning tree case of coefficients
in a flat complex line bundle where $k = dim \ B_{d-1}(X;E_{\rho})$, is discussed
in  M. Catanzaro, V. Chernyak,
and  J. Klein \cite{CatanzaroChernyakKlein1} and there related directly to circuit theory.

\subsection{Manifolds}

When $K$ is a triangulation of a closed n-manifold, one may introduce the dual cell subdivision, $K^\star$
whose $(n-i)$-cells, $D(e_{i,j})$,  are  indexed by the baricenters, say $\hat{\sigma}_{i,j}$ of
the i-dimensional simplices $\sigma_{i,j}$. The dual cell $D(\sigma_{i,j})$ consists of the sets
of simplices  $<\hat{\sigma}_{i,j_1}, \hat{\sigma}_{i+1,j_2}, \cdots, \hat{\sigma}_{n,j_{n-i+1}}>$
of the baricentric subdivision of $K$ where the baricenters correspond to an increasing flag
of simplicies of $K$ starting with $\sigma_{i,j} =\sigma_{i,j_1}$. That is, the simplicies
 $\sigma_{i+k-1,j_{k}}$ of $K$ are of dimension $i-1+k$, with $  k=1, \cdots, n-i+1$ and
 $\sigma_{i,j} = \sigma_{i,j_1} \subset \sigma_{i+1,j_2} \subset \cdots \subset \sigma_{n, j_{n-i+1}}$.
 The bijection $\sigma_{i,j} \mapsto D(\sigma_{i,j})$ establishes an isomorphism  of chains
  $$
 C_i(K;Z) \cong C_{n-i}(K;Z) \ and \ C_i(K;R) \cong C_{n-i}(K;R)
 $$
 under which the boundary $\partial$ of $K$ and its dual $\partial^\star$, the ``coboundary''  for the cellular inner product
 on $K$, are identified with the coboundary for $K^\star$ and the boundary for $K^\star$, respectively.

 Hence, in  this special context, forests of $K$ correspond to dual forests of $K^\star$, coforests of $K$ correspond
 to dual coforests of $K^\star$, etc.

\section{\textbf{Enumerative combinatorics of the simplex  following   Kalai \cite{Kalai}}} \label{sectionKalai}
\vspace{.1in}

Let $\Delta_{n-1}$ denote the standard $n-1$ simplex with $n$ vertices, say $x_1,\cdots, x_n$
in this order. Regard $\Delta_{n-1}$ as  the join
$$
\Delta_{n-1} = x_1 \star \Delta_{n-2}
$$
where $\Delta_{n-2}$ is the standard $n-2$  simplex with vertices $x_2, \cdots,  x_n$.
For a simplex, say $s$ in $\Delta_{n-1}$, let $s$ have its inherited orientation from
this choice of ordering; for emphasis, denote the oriented simplex $\vec{s}$.
As $s$ ranges over the $j$-simplices of $\Delta_{n-1}$, respectively $\Delta_{n-2}$,
the elements $ 1 \cdot \vec{s}$ provide an integral basis for the $j$-chains, $C_j(\Delta_{n-1};Z),
C_j(\Delta_{n-2};Z)$. With respect to these cellular choices of  basis, the boundary maps
$$
\partial_{\Delta_{n-1}} : C_k( \Delta_{n-1};Z) \rightarrow C_{k-1}( \Delta_{n-1};Z)
\ and \ \partial_{\Delta_{n-2}}  : C_k( \Delta_{n-2};Z) \rightarrow C_{k-1}( \Delta_{n-2};Z)
 $$
 are represented by the respective incidence matrices, and have adjoints,
 with respect to the cellular bases,
  $\partial^\star_{\Delta_{n-1}}$ and $\partial_{\Delta_{n-2}}^\star$ represented by the
  transposes of these incidence matrices respectively.

\vspace{.3in}
Note that the integral d-chains, $C_d(\Delta_{n-1};Z)$,  for $d\le n-1$ have a spanning
tree consisting of all $d$-cells, say $s$, containing the initial vertex $x_1$. Moreover,
associated to this spanning tree the integral d-cycles, $Z_d(\Delta_{n-1};Z)$ for $d < n-1$
have a basis consisting of the boundaries $\partial ( x_1 \star t)$ where $t$ ranges over
the d-cells in $\Delta_{n-2}$, that  is the d-cells not containing $x_1$.
Consequently, the above theorems apply in this context.

\vspace{.3in}
 Now consider   the standard cellular inclusion :
 $$
 i : C_{k-1}(\Delta_{n-2};Z) \rightarrow C_{k-1}( \Delta_{n-1};Z)
 $$
In the cellular basis, the adjoint is just the cellular restriction mapping
$$
 C_{k-1}(\Delta_{n-1};Z) \stackrel{ restriction \ = \ i^\star}{\rightarrow } C_{k-1}(\Delta_{n-2};Z) \ .
 $$

 Following, Kalai, consider the composite :
 $$
  i^\star \cdot \partial_{\Delta_{n-1}} : C_k(\Delta_{n-1};Z) \rightarrow C_{k-1}(\Delta_{n-2};Z)
  $$
  It is represented by the restricted incidence matrix, say $I^k_r$,  obtained by deleting  the rows corresponding
  to $k-1$ simplices not in $\Delta_{n-2}$.
  In these terms, one of Kalai's results takes the following form:

  \begin{thm} [Kalai \label{thmKalai1}] The composite
  $$
  (i^\star \cdot \partial_{\Delta_{n-1}}) \cdot ( i^\star \cdot \partial_{\Delta_{n-1}})^\star :
  C_{k-1}( \Delta_{n-2};Z) \rightarrow C_{k-1}( \Delta_{n-2};Z)\ ,
  $$
  which is represented by $I^k_r \cdot (I^k_r)^t$,
  has eigenvalues $1$ and $n$ with respective multiplicities
  $ {n-2}\choose{k-1}$ and ${n-2}\choose{k}$.
  \end{thm}

    Note that, as expected, the sum of these multiplicities is ${n-1}\choose{k}$  = $ rank \ C_{k-1}(\Delta_{n-2};Z)$.

  \vspace{.3in}
  The result easily gives the eigenvalues and their multiplicities of the
  Kirchhoff matrix representing the combinatorial Laplacian
  $$
  \partial_{\Delta_{n-2}} \cdot \partial_{\Delta_{n-2}}^\star  : C_{k-1}(\Delta_{n-2}:Z) \rightarrow   C_{k-1}(\Delta_{n-2}:Z)
  $$
  as follows:

  For $s$ a $k-1$ simplex in $\Delta_{n-2}$, one computes :
  $$
  \begin{array}{l}
  (i^\star \cdot \partial_{\Delta_{n-1}}) \cdot ( i^\star \cdot \partial_{\Delta_{n-1}})^\star ( 1 \cdot \vec{s})
  =  (i^\star \cdot \partial_{\Delta_{n-1}}) \cdot \partial_{\Delta_{n-1}}^\star \cdot i(  1 \cdot \vec{s})
  = (i^\star \cdot \partial_{\Delta_{n-1}}) \cdot (\partial_{\Delta_{n-1}}^\star( 1 \cdot \vec{s}) )  \\
  = (i^\star \cdot \partial_{\Delta_{n-1}}) \cdot (1 \cdot  x_1 \star \vec{s} + i \cdot( \partial_{\Delta_{n-2}}^\star(1 \cdot \vec{s}) ))
  = i^\star (\ 1 \cdot  \vec{s} - x_1 \star \partial_{\Delta_{n-2}}^\star( 1 \cdot \vec{s}\ )
  +  \partial_{\Delta_{n-2}} ( \partial_{\Delta_{n-2}}^\star( 1 \cdot \vec{s}))\  ) \\
  = 1 \cdot \vec{s} +  \partial_{\Delta_{n-2}} ( \partial_{\Delta_{n-2}}^\star( 1 \cdot \vec{s}\ )
  \end{array}
  $$

  That is, one sees that
  $$
  (i^\star \cdot \partial_{\Delta_{n-1}}) \cdot ( i^\star \cdot \partial_{\Delta_{n-1}} )^\star= Id +  \partial_{\Delta_{n-2}} ( \partial_{\Delta_{n-2}}^\star)
  $$
  yielding the following corollary to Kalai's result.
  Let $Bd^k$ be the matrix representing the boundary mapping
  $\partial_{\Delta_{n-2}} :  C_k(\Delta_{n-2}:Z) \rightarrow  C_{k-1}(\Delta_{n-2}:Z)$.

  \begin{corollary} \label{corollaryKalai1} The combinatorial Laplacian
  $$
  \partial_{\Delta_{n-2}} \cdot \partial_{\Delta_{n-2}}^\star  : C_{k-1}(\Delta_{n-2}:Z) \rightarrow   C_{k-1}(\Delta_{n-2}:Z)
  $$
  is represented by the matrix $Bd^{k} \cdot (Bd^{k})^t$ and
 has eigenvalues  $0$ and $n-1$ with respective multiplicities
  $ {n-2}\choose{k-1}$ and ${n-2}\choose{k}$.
 \end{corollary}

 \vspace{.1in}

 Now there is a natural integral and geometric  basis for the integral $k-1$ cycles of $\Delta_{n-1}$, $Z_{k-1}(\Delta_{n-1};Z)$,
 provided  by the following map from $k-1$ chains on $\Delta_{n-2}$ to $k-1$ cycles which is an isomorphism:
 $$
 \begin{array}{l}
 \Phi : C_{k-1}( \Delta_{n-2};Z) \rightarrow Z_{k-1}( \Delta_{n-1};Z) \\
 1 \cdot \vec{s} \mapsto \partial_{\Delta_{n-1}} \ ( x_1 \star \vec{s}) \ \ \ for \ s \ \ a \ \ k \ simplex \ in \ \Delta_{n-2}
 \end{array}
 $$

 Since $\partial_{\Delta_{n-1}} \ ( x_1 \star \vec{s}) = 1 \cdot \vec{s} - x_1 \star \partial_{\Delta_{n-2}}(1 \cdot \vec{s})$,
 the composite $i^\star \Phi ( 1 \cdot \vec{s}) $ equal just $1 \cdot \vec{s}$. That is,
 $
 i^\star \cdot \Phi = Identity.
 $
 Hence, $\Phi$ is one to one.

 Next given a $k-1$ cycle, $z \in Z_{k-1}( \Delta_{n-1};Z)$, then
 $ z= \Phi( i^\star(z) ) + ( z -\Phi( i^\star(z) )$ where $i^\star ( z -\Phi( i^\star(z) ) =0$.
 Hence, $ ( z -\Phi( i^\star(z)) = x_1 \star c$ for some $k-2$ chain $c$ supported in $\Delta_{n-2}$.
 But then $0 = i^\star( \partial_{\Delta_{n-1}}( z -\Phi( i^\star(z)) )=  i^\star( \ \partial_{\Delta_{n-1}}( x_1 \star c) \ )
 = i^\star( \ c - x_1 \star \partial_{\Delta_{n-2}}(c) \ ) = c$. This proves that $c=0$ and so
 that $z= \Phi( i^\star(z) ) $ when $z$ is an integral $k-1$ cycle.

 In toto, this shows that $\Phi$ is an isomorphism.

 In particular, the cellular basis for $C_{k-1}( \Delta_{n-2},;Z)$
 gives the geometric  basis $\{\Phi(1 \cdot \vec{s}) \}$ for the $k-1$ cycles.
 As expected, using the Kalai result again, one may use this method
 to find the eigenvalues and multiplicities of the mesh matrix for
 $k-1$ cycles.

 Let $J^{k-1}$ represent the inclusion $Z_{k-1}(\Delta_{n-1};Z) \subset C_{k-1}(\Delta_{n-1};Z)$
 with respect to the geometrical basis $\{\Phi(s)\}$ of $k-1$ cycles and the cellular
 basis of $k-1$ chains. That is,  $J^{k-1}$  represents the composite $C_{k-1}(\Delta_{n-2};Z) \stackrel{\Phi}{\rightarrow}
 Z_{k-1}(\Delta_{n-1};Z) \stackrel{i}{\rightarrow} C_{k-1}(\Delta_{n-1};Z)$ with respect
 to the cellular bases.

 \begin{corollary} \label{corollaryKalai2} For the integral  basis of $Z_{k-1}(\Delta_{n-1},Z)$ provided by
 the geometric integral basis $\Phi( 1 \cdot \vec{s})$ where $s$ ranges over the $k-1$
 simplices of $\Delta_{n-2}$, the mesh matrix for cycles
 representing  the pairing
 $$
 Z_{k-1}(\Delta_{n-1};Z) \times Z_{k-1}(\Delta_{n-1};Z) \rightarrow Z
 $$
 with respect to  this  basis, or equivalently of the matrix $J^{k-1} (J^{k-1})^t$,
 has eigenvalues only equal to $n$ and $1$ with respective
 multiplicities ${ n-2}\choose{k-1}$ and ${n-2}\choose{k}$.
 \end{corollary}
Since $Z_d(\Delta_n:Z) = B_d(\Delta_n;Z)$, the cycle and boundary mesh matrices
are identical in this case.

 To prove this it suffices to compute and simplify the entries with $s,s^\#$ two
 $k-1$ simplicies in $\Delta_{n-2}$ :
 $$
 \begin{array}{l}
 < \Phi(  1 \cdot \vec{s}) , \Phi( 1 \cdot \vec{s}^\#>
 = < \partial_{\Delta}( x_1 \star \vec{s})  ,  \partial_{\Delta}( x_1 \star \vec{s}^\#) > \\
 = < 1 \cdot \vec{s} - x_1 \star \partial_{\Delta_{n-2}}(\vec{s}) \ , \ 1 \cdot \vec{s}^\#
  - x_1 \star \partial_{\Delta_{n-2}}(\vec{s}^\#) >
  = \delta_{s,s^\#} + <\partial_{\Delta_{n-2}}(\vec{s}), \partial_{\Delta_{n-2}}(\vec{s}^\#)>\\
  = \delta_{s,s^\#} + <\vec{s}, \partial_{\Delta_{n-2}}^\star \partial_{\Delta_{n-2}}(\vec{s}^\#)>
  \end{array}
  $$

Now by Lemma \ref{lemmaalg2} the non-zero eigenvalues of
$ \partial_{\Delta_{n-2}}^\star \partial_{\Delta_{n-2}}  $
acting on $C_{k-1}(\Delta_{n-2};Z)$ are the same as
that of the combinatorial Laplacian, or Kirchhoff map,
$\partial_{\Delta_{n-2}} \partial_{\Delta_{n-2}}^\star$
acting on $C_{k-2}(\Delta_{n-2};Z)$. By corollary
 \ref{corollaryKalai1}, this last has a single non vanishing
 eigenvalue, nameby $n-1$,  of multiplicity $ {n-2}\choose{k-1}$. Hence,
 the mesh matrix above has an  eigenvalue
  $n$ of multiplicity $ {n-2}\choose{k-1}$
 and  eigenvalue $1$  of multiplicity
 $
 {{n-1}\choose{k}} -  {{n-2}\choose{k-1} } = { {n-2}\choose{k} }
 $
 as claimed.

\vspace{.3in}
Kalai \cite{Kalai} has given an additional refinement of
his  result recorded here in  theorem  \ref{thmKalai1}. For each  vertex  $x_j$
introduce a real positive   indeterminant, say $a_j>0$. Introduce a weight for each
face $f=<x_{i_1}, x_{i_2}, \cdots, x_{i_l}>$ with  $w(f)$ equal to the product $a_{i_1} \cdot a_{i_2}\cdots a_{i_l}$
of the corresponding indeterminants of the vertices of the face $f$.
This defines an inner product on each of the $j$ chains
for which the elements $\{ 1 \cdot f /\sqrt{w(f)} \}$ are an orthonormal basis.
These are called the normalized cellular basis,
Let $D(k,m)$ be the diagonal self mapping of the $k$ chains  on $\Delta_m$ which multiplies
a $k$ face $f$ by $w(f)$. For convenience, also let $D(k,m)$ also denote the diagonal matrix which represents
$D(k,m)$ in the cellular basis.

Now  the Kalai mapping $i^\star \cdot \partial_{\Delta_{n-1}}: C_k(\Delta_{n-1};Z) \rightarrow  C_{k-1}(\Delta_{n-2};Z)$
with respect to this  new normalized cellular  basis gives  a reduced incidence matrix, say $(I^k_r)^\#$
which is related to the original reduced incidence matrix $I^k_r$ by
$$
(I^k_r)^\# = D(k-1,n-2)^{-1/2} \ I^k_r \cdot D(k,n-1)^{1/2} \ .
$$

Kalai's extension states that the matrix
$$
(I^k_r)^\# \cdot ((I^k_r)^\#)^t =  D(k-1,n-2)^{-1/2} \ I^k_r \cdot D(k,n-1) \ (I^k_r)^t \cdot D(k,n-2)^{-1/2}
$$
has eigenvalues $a_1 $ and $\Sigma_{j=1}^{n} a_j$ with multiplicities
$ {n-2}\choose{k-1}$ and ${n-2}\choose{k}$ \cite{Kalai}.

\vspace{.3in}
Correspondingly, the boundary mapping $\partial_{\Delta_{n-2}}  : C_{k-1}(\Delta_{n-2}:Z) \rightarrow   C_{k-2}(\Delta_{n-2}:Z)$ with respect to the new normalized cellular basis is a matrix, say $
(Bd^{k-1})^\#$. It is related to the standard boundary mapping $Bd^{k-1}$ by
$$
(Bd^{k-1})^\# = D(k-1,n-2)^{-1/2} \ (Bd^{k-1}) \ D(k, n-1)^{1/2} \ .
$$

The  normalized Kirchhoff mapping is then represented by
$$
(Bd^{k-1})^\# \ ((Bd^{k-1})^\#)^t = D(k-1,n-2)^{-1/2} \ (Bd^{k-1}) \ D(k, n-1) (Bd^{k-1})^t \ D(k-1, n-2)^{-1/2} \ .
$$
A slight extension of the above argument shows that this matrix has eigenvalues $0$ and $\Sigma_{j=1}^{n-1}
\ a_j$ with multiplicities  with respective multiplicities
  $ {n-2}\choose{k-1}$ and ${n-2}\choose{k}$. This follows from Kalai's results by the
  method used in the unnormalized case above.

\vspace{.3in}
Correspondingly, under the identification
$\Phi : C_{k-1}( \Delta_{n-2};Z) \rightarrow Z_{k-1}( \Delta_{n-1};Z)$, one has
a bilinear ``mesh'' pairing
$$
\begin{array}{l}
C_{k-1}( \Delta_{n-2};Z) \times C_{k-1}( \Delta_{n-2};Z) \rightarrow Z_{k-1}(\Delta_{n-1};Z) \times Z_{k-1}(\Delta_{n-1};Z
\rightarrow R
\end{array}
$$
where the last bilinear pairing used the chosen weights $w(f)$.

For the choice of normalized cellular  basis for $C_{k-1} ( \Delta_{n-2};Z)$,
it is easily seen, by the method  above, that this matrix has eigenvalues
$a_1$ and $\Sigma_{j=1}^{n-1} a_j$ of respective multiplicities
  ${ n-2}\choose{k-1}$ and ${n-2}\choose{k}$. This follows from Kalai's results by the
  method used in the unnormalized case above.

\vspace{.1in}
            Extending in several directions some of the formulae of Kalai, the complete calculations in this section
of the eigenvalues of the Kirchhoff matrix and of the two mesh matrices in each dimension of simplices, yield, by
results developed in the preceding sections, many combinatorial identities involving appropriate subobjects
of simplices for each of the coefficients of the corresponding characteristic polynomials.

Sylvain E. Cappell, Courant Institute of Mathematical Sciences, New York University. cappell@cims.nyu.edu

\vspace{.2in}

Edward Y. Miller, Mathematics Department, Tandon School of Engineering of
New York University. em1613@nyu.edu

\end{document}